\documentclass [journal] {IEEEtran}

\usepackage{bm}
\usepackage[english]{babel}
\usepackage{amsmath}
\usepackage[pdftex]{graphicx}

\newcommand{\fraca}[2]{\displaystyle\frac{#1}{#2}}

\newcommand{\mybox}[2]{\makebox[#1\linewidth]{#2}}
\newcommand{\myboxl}[2]{\makebox[#1\linewidth][l]{#2}}

\newcommand{\ap}[2]{\par\noindent
\mybox{0.9}{%
\mybox{0.90}{%
#1\hfill\
\myboxl{0.75}{#2\hfill\mybox{0.05}{\hfill}}}}}

\begin{document}

\title{Inference for a Step-Stress Model
 With Type-II and Progressive Type-II Censoring
 and Lognormally Distributed Lifetimes}

\author{Aida Calvi\~no\thanks{A. Calvi\~no is with the Department
of Statistics and Data Science, Complutense University of Madrid, Madrid, 28040 SPAIN e-mail: aida.calvino@ucm.es.}
}

\markboth{}%
{}

\maketitle

\begin{abstract}
Accelerated life-testing (ALT) is a very useful technique for examining the reliability of highly reliable products. It allows testing the products at higher than usual stress conditions to induce failures more quickly and economically than under typical conditions.
A special case of ALT are step-stress tests that allow experimenter to increase the stress levels at fixed times. This paper deals with the multiple step step-stress model under the cumulative exposure model with lognormally distributed lifetimes in the presence of Type-II and Progressive Type-II censoring.
For this model, the maximum likelihood estimates (MLE) of its parameters, as well as the corresponding observed Fisher Information Matrix (FI), are derived. The likelihood equations do not lead to closed-form expressions for the MLE, and they need to be solved by means of an iterative procedure, such as the Newton-Raphson method. We then evaluate the bias and mean square error of the estimates and provide asymptotic and bootstrap confidence intervals. Finally, in order to asses the performance of the confidence intervals, a Monte Carlo simulation study is conducted.


\end{abstract}

\IEEEpeerreviewmaketitle

\section{Introduction}

\IEEEPARstart{N}{owadays}, most manufactured products are highly reliable with large lifetimes that result in large costs and high experimental times when testing them under typical conditions. In those cases, when conventional life-testing become unuseful, the reliability experimenter may adopt accelerated life-testing, wherein the experimental units are subjected to higher stress levels than under normal operating conditions. Accelerated life tests (ALT) are used to quickly obtain information on the life time distribution of products by testing them at higher than nominal levels of stress to induce early failures (see, for example, \cite{LiuQ:2011} or \cite{XuBT:14}). Furthermore, ALTs allow to examine the effect of stress factors, such as pressure or temperature, on the lifetimes of experimental units. Data collected from ALT needs to be fitted to a model that relates the lifetime to stress and estimate the parameters of the lifetime distribution under normal conditions. This requires a model to relate the levels of stress to the parameters of the lifetime distribution. One such model is the {\em cumulative exposure model} introduced by \cite{Sedyakin:66} and discussed further by \cite{Nelson:80} and \cite{BagdonaviciusN:02}, among others.

Accelerated life-testing may be performed either at increasing or constant high stress levels. In practice, constant stress ALT leads to very few failures within the experimental time, reducing the effectiveness of accelerated testing. A particular case of accelerated testing is the {\em step-stress model}, which allows for a change of stress in steps at various intermediate stages of the experiment\footnote{Stress factors can include humidity, temperature, voltage, load or any other factor that directly affects the life of the products.}. Specifically, a random sample of $n$ units is placed on a life test at an initial stress level $x_1$. At prefixed times $\tau_1,\tau_2,\dots,\tau_{m-1}$, the stress levels are increased to $x_2,x_3,\dots,x_{m}$, respectively. In this paper, we assume both Type-II censoring, where the experiment terminates when a pre-specified number $r$ ($r<n$) of failures is observed (the remaining units are called censored); and Progressive Type-II censoring, where a predetermined number of survival units is removed from the test whenever a failure occurs and, again, the experiment terminates when a pre-specified number of failures is reached.

The step-stress model has been discussed extensively in the literature. \cite{GangulyKM:15}, \cite{Xiong:98}, \cite{BalakrishnanKNK:07}, \cite{BalakrishnanX:07b}, \cite{BalakrishnanX:07a} and \cite{BalakrishnanZX:09} have all considered inferences for the step-stress model assuming exponential lifetimes based on different censoring schemes. \cite{MillerN:83} and \cite{Bai:89} discussed the determination of optimal time $\tau_1$  at which to change the stress level from $x_1$ to $x_2$ in the case of a simple step-stress model ($m=1$). Under progressive Type-II censoring, \cite{XieBD:07} developed both inference and optimal progressive scheme. While all these discussions dealt with exponential step-stress models, \cite{KhamisH:98} and \cite{KateriB:08} examined inferential methods and \cite{NgCB:04} studied optimal progressive plans under Weibull distributed lifetimes. \cite{BalakrishnanZX:09} and \cite{LinC:12} have developed simple and multiple step-stress models, respectively, with lognormally distributed lifetimes
and Type-I censoring. For a comprehensive review on step-stress models refer to \cite{GuonoB:01}, \cite{Nelson:05a} and \cite{Nelson:05b}. \cite{Balakrishnan:09} also provides a detailed review of work on exact inferential procedures for exponential step-stress models as well as associated optimal design of step-stress tests. One may refer to \cite{Balakrishnan:07} for a overview of various developments relating to progressive censoring.

In this paper we deal with the step-stress model with lognormally distributed life times and Type-II and progressive Type-II censoring. In particular, a multiple step-stress model where the location parameters $\mu_i$ of the lognormally distributed lifetimes in each step are given by the linear link function: $\mu_i=\gamma_0+\gamma_1x_i$ is developped. Assuming a link function permits avoiding estimating one location parameter per step, as only $\gamma_0$ and $\gamma_1$ need to be estimated independently of $m$. Furthermore, it allows obtaining the parameters of the lognormal life time distribution for whatever lever of stress $x_i$. \cite{LinC:12} also use a linear link function in their multiple step-stress model, but, as already mentioned, they assume Type-I censoring. When assuming this link function, some physical models, such as the Arrhenius equation or the Inverse Power relationship, can be applied, which permits modeling more real situations. Assumming that the accelerating stress variable $V$ is positive, then the inverse power
relationship is given by:
\[T=AV^\eta,\]
where $A$ and $\eta$ are parameters characteristics of the product, the test method, etc. We can obtain the inverse power law by performing the following transformation: $x_i=\log(V_i)$. According to the Arrhenius rate law (\cite{IUPAC}), the rate of a simple chemical reaction depends on temperature as follow:
\[T=A\exp\left\{\frac{E}{kV}\right\},
\]
where $T$ is the nominal life, $E$ is the activation energy of the reaction, $k$ is Boltzmann's constant, $8.6173\times10^{-5}$ electron-volts per kelvin degree, $V$ is the absolute Kelvin temperature and $A$ is a constant that is characteristic to the product failure mechanism and test conditions. The Arrhenius life-temperature relationship can be obtained by performing the following transformation on our log-link function: $x_i=\frac{1}{kV_i}$. In the examples shown in this paper, we assume the stress level is the temperature and the way it affects the lifetime follows the Arrhenius law.

As already stated, we further assume that the lifetime of a test unit at stress level $x_i$ follows a lognormal distribution, such that:
\[\begin{split}
&\log(T_i)=\gamma_0+\gamma_1x_i+\epsilon,\quad\epsilon\sim N(0,\sigma^2)\\
&E[\log(T_i)]=\gamma_0+\gamma_1x_i=\mu_i\\
&\log(T_i)\sim N(\mu_i,\sigma^2).
\end{split}
\]

The rest of the paper is organized as follows. Section \ref{s2} and \ref{s3} deal with the Multiple Step-Stress model under Type-II and Progressive Type-II censoring, respectively, deriving their associate likelihood function, maximum likelihood estimates and Fisher information matrix. In Section \ref{s4} we show how to construct asymptotic and bootstrap confidence intervals for the parameters. We present an illustrative example in Section \ref{s6}. Section \ref{s5} includes a Monte Carlo simulation study used to evaluate the performance of the methods previously proposed. Finally, in Section \ref{s8} we include some conclusions and future research lines.

\section{Multiple Step-Stress Model under Type-II censoring}\label{s2}

A multiple step-stress testing experiment under a Type-II censoring scheme can be constructed as follows. Let $x_1,x_2,\dots,x_m$ be the stress levels (in the examples below the stress level is temperature), such that $x_1<x_2<\dots<x_m$. A random sample of $n$ experimental units are placed on a life-test at an initial stress level of $x_1$. Then, at a prefixed time $\tau_1$, the stress level is changed to $x_2$; next, at time $\tau_2$, the stress level is changed to $x_3$, and so on. The test is finished when the total number of failed units reaches $r$. For $i=1,2,\dots,m$, let $n_i$ be the number of units failed at stress level $x_i$ (i.e., in the time interval $[\tau_{i-1},\tau_i)$), and $t_{i,j}$ denote the $j$th ordered failure time out of $n_i$ units at level $x_i$, $j=1,2,\dots n_i$.

At stress level $x_i$, the life time $T_i$ of a test unit is assumed to follow a lognormal distribution with location and scale parameter $\mu_i$ and $\sigma$, respectively; whose Cumulative Distribution Function (CDF) and Probability Density Function (PDF) are given by:
\begin{equation}\label{e2}
  F_i(t;\mu_i,\sigma)=\Phi\left(\fraca{\log(t)-\mu_i}{\sigma}\right)
\end{equation}
and
\begin{equation}\label{e3}
  f_i(t;\mu_i,\sigma)=\fraca{1}{\sigma t \sqrt{2\pi}}\exp\left\{-\fraca{1}{2}\left(\fraca{\log(t)-\mu_i}{\sigma}\right)^2\right\},
\end{equation}
respectively.

The location parameters $\mu_i$ are given by the following linear link function
\begin{equation}\label{e2a}
 \mu_i=\mu(x_i)=\gamma_0+\gamma_1 x_i,
\end{equation}
and the scale parameter $\sigma$ is assumed to be free of the stress levels. Therefore, we will need to estimate the regression parameters $\gamma_0$ and $\gamma_1$ (that will give us the relation between the life time and the stress level) as well as the scale parameter.

We further assume that the data comes from a cumulative exposure model, which is the most prominent and commonly used model in the analysis of data observed from step-stress experiments. This model relates the lifetime distribution of experimental units at one stress level to the distributions at preceding stress levels by assuming specifically that the residual life of the experimental units depends only on the cumulative exposure the units have experienced, with no memory of how this exposure was accumulated. The cumulative exposure of the units is measured as the probability of failing up to that time considering the previous stress levels.

At each step $i$, we assume that the lifetime at that specific step (i.e., $t-\tau_{i-1}$, as $\tau_{i-1}$ units of time are spent in the previous steps) follows a log-normal distribution with parameters $\mu_i$ and $\sigma$. However, the units placed in the test have suffered some ``damage'' before step $i$ ( because of the previous exposure to high levels of stress), that needs to be taken into account. This is obtained by adding an artificial extra time ($s_{i-1}$) to the lifetimes at each stress level ($i$) to reflect the exposure suffered at previous levels. As the exposure is reflected in the survival probability, $s_{i-1}$ is given by the equation

\begin{equation}\label{e2d2}
\begin{split}
F_i(s_{i-1};\mu_i;\sigma)=F_{i-1}(\tau_{i-1}+s_{i-2}-&\tau_{i-2};\mu_{i-1},\sigma),\\&i=2,3,\dots,m,
\end{split}
\end{equation}
where $\tau_0=0$, $s_0=0$, and $F_i(t;\mu_i,\sigma)$ is given by Equation (\ref{e2}).

\begin{figure*}
\begin{center}
\includegraphics*[scale=0.6]{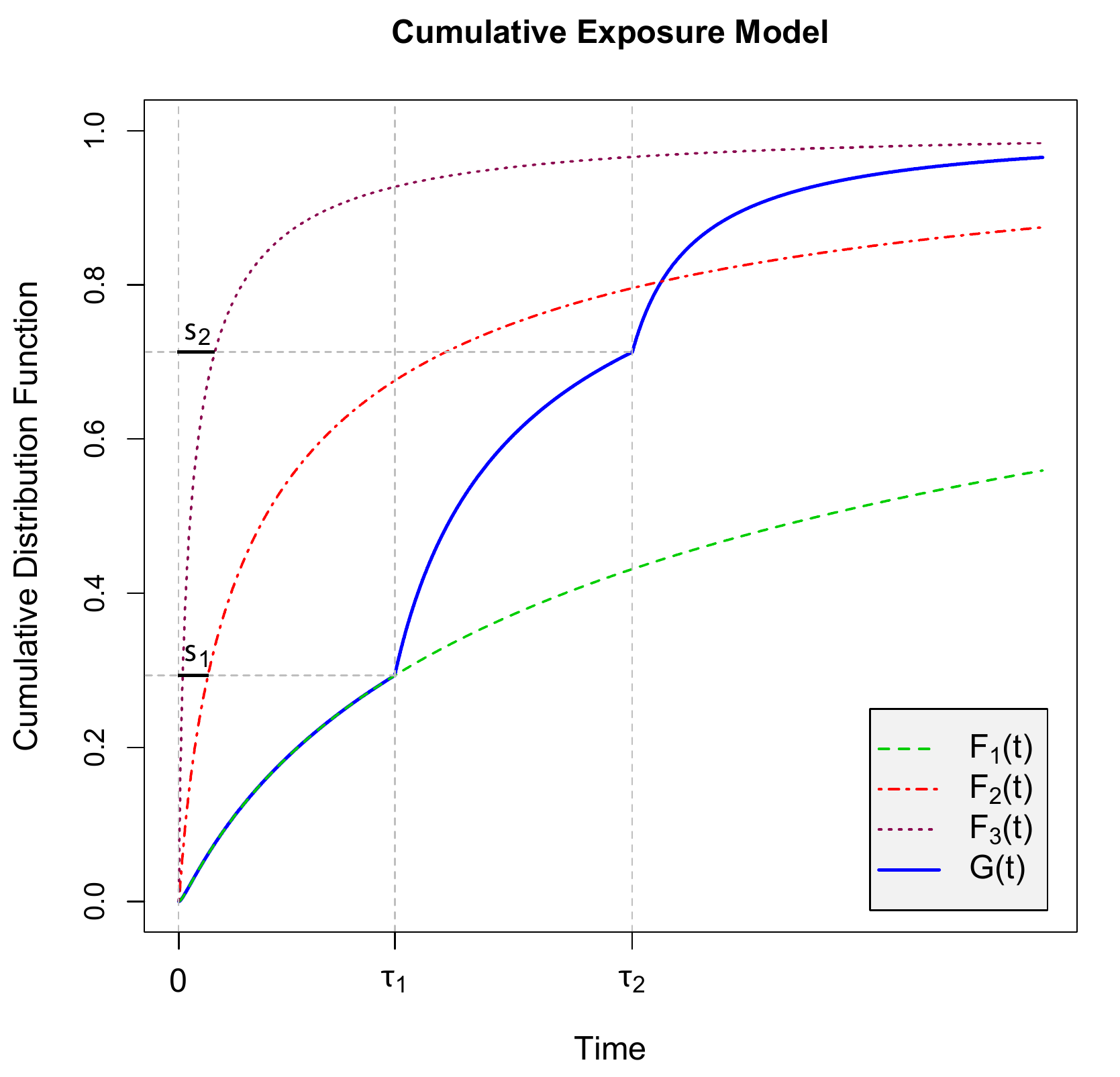}
\caption{\label{fig1}CDF of a test unit under the cumulative exposure model, given in (\ref{e2b}).}
\end{center}
\end{figure*}

Figure \ref{fig1} shows the CDF of a test unit under the cumulative exposure model of a multiple step-stress model ($G(t)$). The construction of the global CDF from the CDF on each step ($F_i(t;\mu_i,\sigma)$) is also shown. Note that the lifetime distributions are shifted in order to get a continuous final function, where the shift parameter is $\tau_{i-1}-s_{i-1}$ for the $i$-th function.

Under the cumulative exposure model, the CDF of the lifetime of a test unit under a multiple step-stress model is given by

\begin{equation}\label{e2b}
\begin{split}
 G(t)=F_i(t+&s_{i-1}-\tau_{i-1};\mu_i,\sigma),\\
&\text{for}\left\{\begin{aligned} &\tau_{i-1}< t\le \tau_i,\quad  i=1,2,\dots,m-1,\\&\tau_{m-1}< t < \infty,\quad  i=m,\end{aligned} \right.
\end{split}
\end{equation}
where $s_{i-1}$ is given by the solution to (\ref{e2d2}), i.e., $s_{i-1}=\left(\tau_{i-1}+s_{i-2}-\tau_{i-2}\right)\exp\left(\mu_i-\mu_{i-1}\right)$ (or equivalently $s_{i-1}=\sum\limits_{h=2}^i(\tau_{h-1}-\tau_{h-2})\exp\{\mu_i-\mu_{h-1}\}$). The PDF of a test unit $g(t)$ can be derived taking derivatives on the corresponding CDF function $G(t)$ in (\ref{e2b}).

\subsection{Likelihood estimation of model parameters}\label{s22}
The likelihood function for the parameter vector $\bm\theta=(\gamma_0,\gamma_1,\sigma)$ based on the observed data is
\begin{equation}\label{e2c}
L(\bm\theta)=\fraca{n!}{(n-r)!}\left\{\prod_{i=1}^{m}\left[\prod_{j=1}^{n_i}g(t_{i,j})\right]\right\}\left[1-G(t_{m,n_m})\right]^{n-r}.
\end{equation}

 Note that there is no need to impose any restriction to the number of failures at each step for the MLE of $\gamma_0,\gamma_1$ and $\sigma$ to exist, even though during some periods (but not more than $m-2$) no failures are observed\footnote{This is an important advantage of using a link function. If we were to estimate the location parameters $\mu_i$ directly, then we will have to impose $n_i>0,\:\forall i$ in order to get all the estimates.}. By substituting (\ref{e2a}) and (\ref{e2b}) in (\ref{e2c}), the log-likelihood function of $\bm\theta$ can be written as

\begin{equation}\label{e2d}
\begin{aligned}
        l(\bm \theta) \propto & -r \log\sigma-\fraca{1}{2}\sum\limits_{i=1}^{m}\sum\limits_{j=1}^{n_i}z_{i,j}^2-\sum\limits_{i=1}^{m}\sum\limits_{j=1}^{n_i}t_{i,j}'\\& +(n-r)\log\left\{1-\Phi\left(z_{m,n_m}\right)\right\},
       \end{aligned}
\end{equation}
where $\log$ is the natural logarithm and
\[\begin{aligned}
t_{i,j}'&=\log\left(t_{i,j}+s_{i-1}-\tau_{i-1}\right),\\
z_{i,j}&=\fraca{t_{i,j}'-\left(\gamma_0+\gamma_1 x_i\right)}{\sigma}.
\end{aligned}\]

Upon differentiating the log-likelihood function of $\bm\theta$ in (\ref{e2d}) with respect to $\gamma_0,\gamma_1$ and $\sigma$, we obtain the following likelihood equations:

\begin{eqnarray}
 \fraca{\partial l(\bm \theta)}{\partial \gamma_0} &=& \fraca{1}{\sigma}\left[\sum\limits_{i=1}^{m}\sum\limits_{j=1}^{n_i}z_{i,j}+(n-r)\fraca{\phi\left(z_{m,n_m}\right)}{1-\Phi\left(z_{m,n_m}\right)}\right] \label{e2e}\\
  \nonumber\fraca{\partial l(\bm \theta)}{\partial \gamma_1} &=& -\sum\limits_{i=1}^{m}\sum\limits_{j=1}^{n_i}z_{i,j}\beta_{i,j}^1-\sum\limits_{i=1}^{m}\sum\limits_{j=1}^{n_i}\beta_{i,j}^2 \\&&-(n-r)\fraca{\phi\left(z_{m,n_m}\right)\beta_{m,n_m}^1}{1-\Phi\left(z_{m,n_m}\right)}\label{e2f}
  \end{eqnarray}
\begin{eqnarray}  
  \nonumber\fraca{\partial l(\bm \theta)}{\partial \sigma} &=& \fraca{1}{\sigma}\left[-r+\sum\limits_{i=1}^{m}\sum\limits_{j=1}^{n_i}z_{i,j}^2\right.\\&&+\left.(n-r)\fraca{z_{m,n_m}\phi\left(z_{m,n_m}\right)}{1-\Phi\left(z_{m,n_m}\right)}\right]\label{e2g}
\end{eqnarray}
where
\begin{eqnarray*}
 \beta_{i,j}^1 &=& \fraca{\partial z_{i,j}}{\partial \gamma_1}=\fraca{\beta_{i,j}^2-x_i}{\sigma},\\
 \beta_{i,j}^2 &=& \fraca{\partial t_{i,j}'}{\partial \gamma_1}=\fraca{\beta_i^3}{t_{i,j}+s_{i-1}-\tau_{i-1}},\\
 \beta_{i}^3 &=& \fraca{\partial s_{i-1}}{\partial \gamma_1}=\sum\limits_{h=2}^i\left(\tau_{h-1}-\tau_{h-2}\right)\left(x_i-x_{h-1}\right)\\&&\times\exp\{\gamma_1(x_i-x_{h-1})\}.
 \end{eqnarray*}

The MLEs of $\bm \theta=(\gamma_0,\gamma_1,\sigma)$ are then obtained by simultaneously solving the following nonlinear equations:
\begin{equation}\label{e13}
  \begin{pmatrix}
    \fraca{\partial l(\bm \theta)}{\partial \gamma_0}\\\fraca{\partial l(\bm \theta)}{\partial \gamma_1}\\\fraca{\partial l(\bm \theta)}{\partial \sigma}
  \end{pmatrix}
  =\begin{pmatrix}
    0\\0\\0
  \end{pmatrix}.
\end{equation}

Since explicit solutions do not exist for (\ref{e2e})-(\ref{e2g}), numerical methods such as the Newton-Raphson procedure need to be used to compute the MLE of $\bm \theta$.

\subsection{Fisher Information Matrix}\label{s23}

The Fisher information matrix, which is the inverse of the variance-covariance matrix of the MLE of the vector parameter $\bm \theta$, denoted by ${\cal F}(\bm \theta)$, is given by
\begin{equation}\label{e14}
  {\cal F}(\bm \theta)=E\left\{\left[\fraca{\partial}{\partial \bm \theta}l(\bm \theta)\right]^2\right\}=-E\left[\fraca{\partial^2}{\partial \bm \theta^2}l(\bm \theta)\right],
\end{equation}
and can be approximated by the observed fisher information matrix ${\cal I}_{obs}(\bm \theta)$ which is the symmetric matrix
\begin{equation}\label{e14a}
  {\cal I}_{obs}(\bm \theta)={\cal I}_{obs}(\gamma_0,\gamma_1,\sigma)=(O_{ij})=\left(-\fraca{\partial^2}{\partial\theta_i\partial\theta_j}l(\bm\theta)\right),
\end{equation}
where $\bm\theta=(\theta_1,\theta_2,\theta_3)=(\gamma_0,\gamma_1,\sigma)$. Given the log-likelihood function $l(\bm\theta)$ in (\ref{e2d}), the elements of the observed fisher information matrix can be found in (\ref{eqA})-(\ref{eqB}) at the top of the following page, where
\begin{eqnarray*}
\beta_{i,j}^4 &=& \fraca{\partial \beta_{i,j}^1}{\partial \gamma_1}= \fraca{\beta_{i,j}^5}{\sigma}\\
\beta_{i,j}^5 &=& \fraca{\partial \beta_{i,j}^2}{\partial \gamma_1}= \fraca{\left(t_{i,j}+s_{i-1}-\tau_{i-1}\right)}{\left(t_{i,j}+s_{i-1}-\tau_{i-1}\right)^2}\sum\limits_{h=2}^i\left(\tau_{h-1}-\tau_{h-2}\right)\\&&\times(x_i-x_{h-1})^2\exp\{\gamma_1(x_i-x_{h-1})\}-\left(\beta_i^3\right)^2.
 \end{eqnarray*}
 \newcounter{MYtempeqncnt}
 
 \begin{figure*}[!t]
\normalsize
\setcounter{MYtempeqncnt}{\value{equation}}
\setcounter{equation}{16}
\begin{eqnarray}
 O_{11}=\fraca{\partial^2 l(\bm \theta)}{\partial \gamma_0^2} \!\!\!&=&\!\!\!  \fraca{1}{\sigma^2}\left\{-r+(n-r)\fraca{\phi(z_{m,n_m})}{1-\Phi(z_{m,n_m})}\left[z_{m,n_m}-\fraca{\phi(z_{m,n_m})}{1-\Phi(z_{m,n_m})}\right]\right\}\label{eqA} \\
 O_{22}=\fraca{\partial^2 l(\bm \theta)}{\partial \gamma_1^2}\!\!\! &=& \!\!\! -\sum\limits_{i=1}^{m}\sum\limits_{j=1}^{n_i}\left[(\beta_{i,j}^1)^2+z_{i,j}\beta_{i,j}^4+\beta_{i,j}^{5}\right]\nonumber\\&&-(n-r)\fraca{\phi(z_{m,n_m})}{1-\Phi(z_{m,n_m})}\left[\beta_{m,n_m}^4-(\beta_{m,n_m}^1)^2\left(z_{m,n_m}-\fraca{\phi(z_{m,n_m})}{1-\Phi(z_{m,n_m})}\right)\right] \\
 O_{33}=\fraca{\partial^2 l(\bm \theta)}{\partial \sigma^2} \!\!\!&=&\!\!\!-\fraca{1}{\sigma^2}\left\{-r+3\sum\limits_{i=1}^{m}\sum\limits_{j=1}^{n_i}z_{i,j}^2+(n-r)\fraca{z_{m,n_m}\phi(z_{m,n_m})}{1-\Phi(z_{m,n_m})}\left[2-z_{m,n_m}\left(z_{m,n_m}-\fraca{\phi(z_{m,n_m})}{1-\Phi(z_{m,n_m})}\right)\right]\right\}\\
 O_{12}=\fraca{\partial^2 l(\bm \theta)}{\partial \gamma_0\partial\gamma_1}\!\!\! &=& \!\!\! \fraca{1}{\sigma}\left\{\sum\limits_{i=1}^{m}\sum\limits_{j=1}^{n_i}\beta_{i,j}^1-(n-r)\fraca{\beta_{m,n_m}^{1}\phi(z_{m,n_m})}{1-\Phi(z_{m,n_m})}\left[z_{m,n_m}-\fraca{\phi(z_{m,n_m})}{1-\Phi(z_{m,n_m})}\right]\right\} \\
  O_{13}=\fraca{\partial^2 l(\bm \theta)}{\partial\gamma_0\partial \sigma} \!\!\!&=&\!\!\! -\fraca{1}{\sigma^2}\left\{(n-r)\fraca{\phi(z_{m,n_m})}{1-\Phi(z_{m,n_m})}\left[1-z_{m,n_m}\left(z_{m,n_m}-\fraca{\phi(z_{m,n_m})}{1-\Phi(z_{m,n_m})}\right)\right]+2\sum\limits_{i=1}^{m}\sum\limits_{j=1}^{n_i}z_{i,j}\right\}\\
  O_{23}=\fraca{\partial^2 l(\bm \theta)}{\partial\gamma_1\partial\sigma}\!\!\! &=& \!\!\! \fraca{1}{\sigma}\left\{(n-r)\fraca{\beta_{m,n_m}^1\phi(z_{m,n_m})}{1-\Phi(z_{m,n_m})}\left[1-z_{m,n_m}\left(z_{m,n_m}-\fraca{\phi(z_{m,n_m})}{1-\Phi(z_{m,n_m})}\right)\right]		+2\sum\limits_{i=1}^{m}\sum\limits_{j=1}^{n_i}z_{i,j}\beta_{i,j}^1\right\}\label{eqB}
\end{eqnarray}
\setcounter{equation}{\value{MYtempeqncnt}}
\hrulefill
\vspace*{4pt}
\end{figure*}

{\centering\section{Multiple Step-Stress Model under Progressive Type-II Censoring}\label{s3}}

In this section a multiple step-stress model under progressive Type-II censoring is developed. A progressively Type-II censored sample is observed as follows. $n$ identical units are placed on a life-testing experiment, and $r$ and $R_k$ ($k=1,\dots,r$) are fixed in advance, such that $R_r=n-r-R_1-\dots-R_{r-1}$. At the time of the first failure, $R_1$ of the $n-1$ surviving units are randomly removed from the experiment; at the time of the second failure, $R_2$ of the $n-2-R_1$ surviving units are randomly removed from the experiment, and so on; the test continues until the $r^{th}$ failure occurs at which time all the remaining $R_r$ surviving units are removed. If $R_1=\dots=R_r=0$, then $n=r$ which corresponds to the complete sample situation. If $R_1=\dots=R_{r-1}=0$, then $R_r=n-r$ which corresponds to the conventional Type-II censoring scheme utilized in Section \ref{s2}.

The main advantage of Progressive Censoring is that it permits obtaining more information about the lifetimes of the units without the need of exposing all units to high levels of stress, and therefore, with lower associated costs.

Again, at stress level $x_i$, the life time $T_i$ of a test unit is assumed to follow a lognormal distribution with parameters $\mu_i$ and $\sigma$ and CDF and PDF given in (\ref{e2}) and (\ref{e3}), respectively. We further assume the cumulative exposure model and, hence, the CDF of the lifetime of a test unit under a multiple step-stress model is given by (\ref{e2b}).

Following the lines of Sections \ref{s22} and \ref{s23}, the likelihood function for the parameter vector $\bm\theta=(\gamma_0,\gamma_1,\sigma)$ based on the observed data is

\begin{equation}\label{e2c2}
L(\bm\theta)=C\prod_{i=1}^{m}\left\{\prod_{j=1}^{n_i}g(t_{i,j})\left[1-G(t_{i,j})\right]^{R_{k(i,j)}}\right\},
\end{equation}
where $k(i,j)$ refers to the ordered position of the $(i,j)$-th unit in the global sample and $C$ is given by:
\[C=\prod_{j=1}^{r}\left\{\sum_{k=j}^{r}(R_k+1)\right\}.\]

Therefore, the maximum likelihood estimates $\hat\gamma_0,\hat\gamma_1$ and $\hat\sigma$ can be obtained by solving the following likelihood equations:

\begin{eqnarray}
 \nonumber\fraca{\partial l(\bm \theta)}{\partial \gamma_0}\!\!\!\! &=& \!\!\!\! \fraca{1}{\sigma}\left[\sum\limits_{i=1}^{m}\sum\limits_{j=1}^{n_i}z_{i,j}+\sum\limits_{i=1}^{m}\sum\limits_{j=1}^{n_i}R_{k(i,j)}\fraca{\phi\left(z_{i,j}\right)}{1-\Phi\left(z_{i,j}\right)}\right]=0 \\
 && \label{e4c}\\
 \nonumber \fraca{\partial l(\bm \theta)}{\partial \gamma_1}\!\!\!\! &=& \!\!\!\! -\sum\limits_{i=1}^{m}\sum\limits_{j=1}^{n_i}z_{i,j}\beta_{i,j}^1-\sum\limits_{i=1}^{m}\sum\limits_{j=1}^{n_i}\beta_{i,j}^2 \\ \!\!\!\!&&\!\!\!\!-\sum\limits_{i=1}^{m}\sum\limits_{j=1}^{n_i}R_{k(i,j)}\fraca{\phi\left(z_{i,j}\right)\beta_{i,j}^1}{1-\Phi\left(z_{i,j}\right)}=0\label{e4d}
 \end{eqnarray}
 \setcounter{equation}{22}
 \begin{eqnarray}
 \nonumber \fraca{\partial l(\bm \theta)}{\partial \sigma}\!\!\!\! &=& \!\!\!\! \fraca{1}{\sigma}\left[-r+\sum\limits_{i=1}^{m}\sum\limits_{j=1}^{n_i}z_{i,j}^2\right.\\&&+\left.\sum\limits_{i=1}^{m}\sum\limits_{j=1}^{n_i}R_{k(i,j)}\fraca{z_{i,j}\phi\left(z_{i,j}\right)}{1-\Phi\left(z_{i,j}\right)}\right]=0,\label{e4e}
\end{eqnarray}
where $t_{i,j}'$, $z_{i,j}$, $\beta_{i,j}^1$ and $\beta_{i,j}^2$ are as given in the preceding section.

Again, as explicit solutions do not exist for (\ref{e4c})-(\ref{e4e}), numerical methods such as the Newton-Raphson procedure need to be used to compute the MLE of $\bm\theta$.

Finally, the elements of the observed fisher information matrix ${\cal I}_{obs}(\bm \theta)=(O_{ij})$ as defined in (\ref{e14a}), where $\bm\theta=(\theta_1,\theta_2,\theta_3)=(\gamma_0,\gamma_1,\sigma)$ can be found in (\ref{eqC})-(\ref{eqD}) at the top of the following page, where $\beta_{i,j}^4$ and $\beta_{i,j}^5$ are as given in the preceding section.
 \newcounter{MYtempeqncnt2}

 \begin{figure*}[!t]
\normalsize
\setcounter{MYtempeqncnt2}{\value{equation}}
\setcounter{equation}{25}
\begin{eqnarray}
 \fraca{\partial^2 l(\bm \theta)}{\partial \gamma_0^2} &=&  \fraca{1}{\sigma^2}\left\{-r+\sum\limits_{i=1}^{m}\sum\limits_{j=1}^{n_i}R_{k(i,j)}\fraca{\phi(z_{i,j})}{1-\Phi(z_{i,j})}\left[z_{i,j}-\fraca{\phi(z_{i,j})}{1-\Phi(z_{i,j})}\right]\right\} \label{eqC}\\
 \fraca{\partial^2 l(\bm \theta)}{\partial \gamma_1^2} &=&  -\sum\limits_{i=1}^{m}\sum\limits_{j=1}^{n_i}\left[(\beta_{i,j}^1)^2+z_{i,j}\beta_{i,j}^4+\beta_{i,j}^{5}\right]\nonumber\\&&-\sum\limits_{i=1}^{m}\sum\limits_{j=1}^{n_i}R_{k(i,j)}\fraca{\phi(z_{i,j})}{1-\Phi(z_{i,j})}\left[\beta_{i,j}^4-(\beta_{i,j}^1)^2\left(z_{i,j}-\fraca{\phi(z_{i,j})}{1-\Phi(z_{i,j})}\right)\right] \\
 \fraca{\partial^2 l(\bm \theta)}{\partial \sigma^2} &=&-\fraca{1}{\sigma^2}\left\{-r+3\sum\limits_{i=1}^{m}\sum\limits_{j=1}^{n_i}z_{i,j}^2+\sum\limits_{i=1}^{m}\sum\limits_{j=1}^{n_i}R_{k(i,j)}\fraca{z_{i,j}\phi(z_{i,j})}{1-\Phi(z_{i,j})}\left[2-z_{i,j}\left(z_{i,j}-\fraca{\phi(z_{i,j})}{1-\Phi(z_{i,j})}\right)\right]\right\}\\
 \fraca{\partial^2 l(\bm \theta)}{\partial \gamma_0\partial\gamma_1} &=&  \fraca{1}{\sigma}\left\{\sum\limits_{i=1}^{m}\sum\limits_{j=1}^{n_i}\beta_{i,j}^1-\sum\limits_{i=1}^{m}\sum\limits_{j=1}^{n_i}R_{k(i,j)}\fraca{\beta_{i,j}^{1}\phi(z_{i,j})}{1-\Phi(z_{i,j})}\left[z_{i,j}-\fraca{\phi(z_{i,j})}{1-\Phi(z_{i,j})}\right]\right\} \\
  \fraca{\partial^2 l(\bm \theta)}{\partial\gamma_0\partial \sigma} &=& -\fraca{1}{\sigma^2}\left\{2\sum\limits_{i=1}^{m}\sum\limits_{j=1}^{n_i}z_{i,j}+\sum\limits_{i=1}^{m}\sum\limits_{j=1}^{n_i}R_{k(i,j)}\fraca{\phi(z_{i,j})}{1-\Phi(z_{i,j})}\left[1-z_{i,j}\left(z_{i,j}-\fraca{\phi(z_{i,j})}{1-\Phi(z_{i,j})}\right)\right]\right\}\\
  \fraca{\partial^2 l(\bm \theta)}{\partial\gamma_1\partial\sigma} &=& \fraca{1}{\sigma}\left\{2\sum\limits_{i=1}^{m}\sum\limits_{j=1}^{n_i}z_{i,j}\beta_{i,j}^1+\sum\limits_{i=1}^{m}\sum\limits_{j=1}^{n_i}R_{k(i,j)}\fraca{\beta_{i,j}^1\phi(z_{i,j})}{1-\Phi(z_{i,j})}\left[1-z_{i,j}\left(z_{i,j}-\fraca{\phi(z_{i,j})}{1-\Phi(z_{i,j})}\right)\right]\right\} \label{eqD}
\end{eqnarray}
\setcounter{equation}{\value{MYtempeqncnt2}}
\hrulefill
\vspace*{4pt}
\end{figure*}

\section{Interval Estimation of Model Parameters}\label{s4}

In this section, we propose different methods of constructing confidence intervals (CI) for the parameter vector $\bm \theta$. In particular, approximate and bootstrap CI are shown as the exact distributions of the MLEs are not explicit.

\subsection{Normal approximation of the MLE}
The approximate confidence intervals for the unknown model parameters are based on the asymptotic distribution of the MLEs of $\bm\theta$ given by
\begin{equation*}
  \bm\hat\theta-\bm\theta\xrightarrow{d} N(0,v(\bm\theta))
\end{equation*}
where $v(\bm\theta)$ is the Cramer-Rao lower bound given by the inverse matrix of ${\cal F}(\bm \theta)$. Since an explicit expression can not be derived for ${\cal F}(\bm \theta)$, we use the observed information matrix ${\cal I}_{obs}(\bm \theta)$ to estimate ${\cal F}(\bm \theta)$, as already mentioned. Let
\begin{equation}\label{e3a}
{\cal I}^{-1}_{obs}(\bm \hat\theta)=\begin{pmatrix}
    V_{11} & V_{12} & V_{13}\\V_{21} & V_{22} & V_{23}\\V_{31} & V_{32} & V_{33}
  \end{pmatrix}.\end{equation}
be the variance-covariance matrix of $\bm\hat\theta$, with ${\cal I}_{obs}$ as given in Sections \ref{s2} and \ref{s3}, for the cases of Type-II and progressive Type-II censoring schemes, respectively.

Based on Slutsky's Theorem, we can show that the pivotal quantities $Z_{\theta_i}=(\hat\theta_i-\theta_i)/\sqrt{V_{ii}},\;i=1,2,3$, converge in distribution to the standard normal distribution. Therefore, two-sided $100(1-\alpha)\%$ approximate CIs for $\bm\theta$ are given by:
\begin{equation}
  \hat\theta_i\pm z_{1-\alpha/2}\sqrt{V_{ii}},\; i=1,2,3,
\end{equation}
where $z_{1-\alpha/2}$ is the upper $\alpha/2$ percentile of the standard normal distribution. Note that these CIs are valid for both Type-II and progressive Type-II censoring schemes by changing the observed fisher information matrix ${\cal I}_{obs}(\bm\theta)$.

\subsection{Bootstrap Confidence Intervals}

The bootstrap method, which is one of a broader class of resampling methods, uses Monte Carlo sampling to generate an empirical sampling distribution of the
estimate $\bm\hat\theta$ (see \cite{EfronT:93} for more details on bootstrap methods). Therefore, before explaining how the limits of the confidence intervals are computed, we first present the algorithm used for simulating the required bootstrap Type-II and progressive Type-II censored data.
\subsubsection*{Bootstrap Sample}
\begin{description}
  \item[Step 1. ] Compute the MLE $\bm{\hat\theta}$ of the parameters, given the censored data along with the pre-fixed times $\tau_1,\dots,\tau_{m-1}$.
  \item[Step 2. ] For the case of Type-II censoring, generate $n$ random observations from a Uniform(0,1) distribution, sort them and keep the first $r$ simulations. Contrary, for the case of progressive Type-II censoring, use the algorithm of \cite{BalakrishnanS:95}. Finally, for both cases denote the corresponding order statistics by $\{U_{1:r},\dots,U_{r:r}\}$.
  \item[Step 3. ] Find $n_1$ such that
  \[U_{n_1:r}<\Phi\left(\fraca{\log(\tau_1)-(\hat\gamma_0+\hat\gamma_1x_1)}{\hat\sigma}\right)<U_{n_1+1:r}.
  \]
  For $l=1,\dots,n_1$, set
  \[t^*_{1,l}=\exp\{\hat\sigma\Phi^{-1}(U_{l:r})+(\hat\gamma_0+\hat\gamma_1x_1)\}.\]
  \item[Step 4. ] For $i=2,\dots m-1$, find $n_i$, such that
  \[\begin{split}
 U_{\sum\limits_{j=1}^{i}n_j:r}<\Phi&\left\{\fraca{1}{\hat\sigma}\left(\log(\tau_i+\hat s_{i-1}-\tau_{i-1})\right.\right.\\&\left.\phantom{\fraca{1}{\hat\sigma}}\left.-(\hat\gamma_0+\hat\gamma_1x_1)\right)\right\}<U_{1+\sum\limits_{j=1}^{i}n_j:r}, \end{split}\]
  where $\hat s_{i-1}$ is the MLE of $s_{i-1}$. Furthermore, for $l=\sum\limits_{j=1}^{i-1}n_j+1,\dots,\sum\limits_{j=1}^{i}n_j$ set
  \[\begin{split}
  t^*_{i,l-\sum\limits_{j=1}^{i-1}n_j}=\tau_{i-1}&-\hat s_{i-1}+\exp\left\{\hat\sigma\Phi^{-1}(U_{l:r})\right.\\&+\left.(\hat\gamma_0+\hat\gamma_1x_1)\right\}.
  \end{split}\]
  \item[Step 5. ] Set $n_m=r-\sum\limits_{j=1}^{m-1}n_j$ and for $l=\sum\limits_{j=1}^{m-1}n_j+1,\dots,r$ set
  \[\begin{split}t^*_{m,l-\sum\limits_{j=1}^{m-1}n_j}&=\exp\left\{\hat\sigma\Phi^{-1}(U_{l:r})+(\hat\gamma_0+\hat\gamma_1x_m)\right\}\\&+\tau_{m-1}-\hat s_{m-1}.
  \end{split}\]
  \item[Step 6. ] Compute the MLEs $\hat\gamma_0^*$, $\hat\gamma_1^*$ and $\hat\sigma^*$ based on $(t^*_{1,1},\dots,t^*_{1,n_1},\dots,t^*_{m,n_m})$.
  \item[Step 7. ] Repeat Steps 2-6 $B$ times and arrange all $\hat\gamma_0^*$, $\hat\gamma_1^*$ and $\hat\sigma^*$ values in
ascending order to obtain the corresponding bootstrap sample
\[\{\hat\theta_k^{*[1]},\hat\theta_k^{*[2]},\dots,\hat\theta_k^{*[B]}\},\quad k=1,2,3,
\]
where $\bm{\hat\theta}^*=(\hat\gamma_0^*, \hat\gamma_1^*, \hat\sigma^*)$.
\end{description}

\subsubsection*{Percentile Bootstrap CI} The percentile interval uses the quantiles of the bootstrap distribution to obtain a confidence interval. In particular,
a two-sided $100(1-\alpha)\%$ percentile bootstrap confidence interval for $\theta_k$ is
\[(\hat\theta_k^{*\left[B\frac{\alpha}{2}\right]},\hat\theta_k^{*\left[B\left(1-\frac{\alpha}{2}\right)\right]}),\quad k=1,2,3.\]

\section{Illustrative example}\label{s6}

Let us consider the simulated data presented in Table \ref{t5} for a $3$-step step-stress model with sample size $n=35$, $20\%$ of censored data and $\tau_1=95$, and $\tau_2=97.5$ for the choice of parameters $\gamma_0=0.76, \gamma_1=0.107$ and $\sigma=0.05$. As it will be shown in the following section, this parameters have been selected to represent a real situation where the temperature is the stress factor.

\begin{table*}
 \caption{Simulated lognormal $3$-step step-stress data of size $n=35$ with fixed times $\tau_1=95$, and $\tau_2=97.5$ and parameters $\gamma_0=0.76$, $\gamma_1=0.107$ and $\sigma=0.05$.}
 \begin{center}
 \begin{tabular}{ccccccccc}
 \hline
 Parameters & \multicolumn{8}{c}{Lifetimes}\\
 \hline
 $\mu_1=4.60$ & 89.406 & 92.317 & 92.651 & 93.755 & 94.483 & 94.985 & &\\
 $\mu_2=3.69$ & 95.018 & 95.218 & 95.352 & 95.441 & 95.461 & 95.835 & 95.854 & 95.903\\
              & 96.321 & 96.430 & 96.508 & 96.568 & 97.206 & 97.463 & &\\
 $\mu_3=2.98$ & 97.509 & 97.604 & 97.971 & 98.070 & 98.104 & 98.202 & 98.278 & 98.507\\
              & 98.548 & 98.549 & 98.565 & 98.710 & 98.861 & 98.880 & 99.058 &\\
 \hline
 \end{tabular}
 \end{center}
 \label{t5}
 \end{table*}

For these data, the MLEs of $\gamma_0, \gamma_1$ and $\sigma$ and the corresponding standard errors were all determined from the formulas in Section \ref{s3} for different censoring schemes, where the censored units have been selected randomly among the possible ones\footnote{For convenience, we have used a special notation in the tables for the progressive censoring schemes. For example, $(3^\star0,4)$ denotes the progressive censoring scheme $(0,0,0,4)$.}. The approximate and bootstrap confidence intervals, as well as the estimates and standard errors, are shown in Table \ref{t6}.
As it can be seen in the tables, the smallest standard errors are obtained for the progressive censoring schemes (in particular, they achieve the smallest value in the case where the censorship takes place in the middle values). Moreover, note that the worst results are obtained for the classical Type-II censoring, as the standard errors are larger and the confidence intervals wider.

Finally, a hypothesis test over $\gamma_1$ has been conducted to test if it is significatively positive ($H_0: \gamma_1\le 0;\; H_1:\gamma_1>0$), i.e., if the stress level has a real influence on the lifetime\footnote{Note that we test the positiveness of the parameter as the considered stress level $1/(kV_i)$ decreases over the levels.}.
For that purpose, a t-test using as the variance the corresponding value in the observed FI matrix and a bootstrap test based on unilateral CI were carried out leading to p-values smaller than 0.01 in both cases. Therefore, as expected, the parameter $\gamma_1$ is significatively positive or, in other words, temperature influences lifetime of test units.

In order to evaluate the convenience of an accelerated-life test, it is important to test the effect of the stress to confirm that it is the cause of the earlier failures.

\begin{table*}
 \caption{Estimated values and confidence intervals of parameters $\gamma_0$, $\gamma_1$ and $\sigma$ for the data in Table \ref{t5} under different censoring schemes.}
 \begin{center}
    \resizebox{!}{2.5cm}{
 \renewcommand{\arraystretch}{1.2}
 \renewcommand{\tabcolsep}{0.9mm}
 \begin{tabular}{cccccccccc}
 \hline
 \multicolumn{4}{c}{} & \multicolumn{2}{c}{$90\%$ CI} & \multicolumn{2}{c}{$95\%$ CI} & \multicolumn{2}{c}{$99\%$ CI}\\
 \hline
 Cens. scheme &  & $\hat\theta_i$ & $\hat{SE}(\hat\theta_i)$ & App. & Boot. & App. & Boot. & App. & Boot.\\
\hline
 $(7,27^\star0)$ & $\gamma_0$ & 0.796 & 1.174&(-1.135, 2.728) & (-1.074, 2.978) & (-1.505, 3.098) & (-1.417, 3.610) & (-2.228, 3.821) & (-1.956, 4.017)\\
& $\gamma_1$  & 0.106 & 0.033 & (0.052, 0.161) & (0.043, 0.159) &(0.041, 0.171)& (0.027, 0.169) & (0.021, 0.192)&(0.016, 0.185)\\
& $\sigma$  & 0.050 & 0.018 & (0.021, 0.080) & (0.023, 0.081) &(0.015, 0.085)& (0.019, 0.088) &(0.004, 0.096) & (0.015, 0.11)\\
\hline
$(27^\star0,7)$ & $\gamma_0$  & 0.270 & 1.270 & (-1.819, 2.358)& (-4.665, 0.556) & (-2.29, 2.758) & (-5.572, 1.002)& (-3.001, 3.540) &(-7.295, 1.940)\\
& $\gamma_1$  &0.121 & 0.036 & (0.062, 0.180) & (0.112, 0.262) &(0.051, 0.191) &(0.100, 0.287) & (0.028, 0.213)& (0.074, 0.339)\\
& $\sigma$  & 0.054 & 0.018 & (0.024, 0.083) & (0.039, 0.129) &(0.019, 0.089) &(0.036, 0.153) &(0.008, 0.100) & (0.027, 0.176)\\
\hline
$(7^\star(0,0,1,0))$ & $\gamma_0$  &1.084& 1.129 &(-0.770, 2.939)&(-1.074, 3.140) &(-1.125, 3.294) &(-1.500, 3.648) &(-1.820, 3.989) &(-2.324, 4.456)\\
& $\gamma_1$  & 0.098 & 0.032 & (0.046, 0.151) & (0.040, 0.157)&(0.036, 0.161) &(0.026, 0.172) &(0.016, 0.180) &(0.003, 0.195)\\
& $\sigma$  & 0.050 & 0.016 & (0.024, 0.077) & (0.025, 0.082)&(0.019, 0.082) &(0.022, 0.091) &(0.009, 0.092) &(0.014, 0.098)\\
\hline
$(10^\star0,7^\star1,11^\star0)$ & $\gamma_0$  & 1.276 & 1.120 &(-0.566, 3.118) &(-0.426, 3.614) &(-0.918, 3.471) &(-0.702, 4.158) &(-1.608, 4.161) & (-1.041, 4.992)\\
& $\gamma_1$  & 0.093 & 0.032 &(0.041, 0.145) &(0.026, 0.139) &(0.031, 0.155) & (0.012, 0.149)&(0.011, 0.174) & (-0.012, 0.158)\\
& $\sigma$  & 0.050 & 0.016 & (0.023, 0.077) & (0.024, 0.075)&(0.018, 0.082) & (0.021, 0.083) & (0.008, 0.092) & (0.014, 0.093)\\
 \hline
 \end{tabular}}
 \end{center}
 \label{t6}
 \end{table*}

\section{Simulation study}\label{s5}

In order to evaluate the performance of the methods described in the preceding sections, i.e., point and interval estimation, a Monte Carlo simulation study has been conducted.
In this simulation study, we have used two sample sizes $n=35$ and $75$ with different proportions of censoring, e.g. $60\%$ and $20\%$, and three censoring schemes (corresponding to: a) complete left censored; b) classical Type-II or complete right censored; and c) balanced progressive scheme, where we have used the special notation explained in Section \ref{s6}) and two number of steps: $m=2$ and $m=3$.

So as to represent a real situation, we have used the parameters of the Arrhenius law. In particular, we have assumed the temperature stress levels (in Celsius degrees) and the mean lifetimes (in units of time) to be: $x_1=50^{\circ}{\rm C},\;E[T_1]=100$, $x_2=150^{\circ}{\rm C},\;E[T_2]=40$ and $x_3=300^{\circ}{\rm C},\;E[T_3]=20$.
So as to evaluate the influence of large variance, two different values of $\sigma$ (0.05 and 0.2) are considered, which are derived from the assumption that the standard deviation of the lifetime at the first stress level are $5$ and $20$ units of time, respectively. Moreover, the values of $\gamma_0=0.76$ and $0.71$ and $\gamma_1=0.107$ and $0.108$ are obtained by applying the formulas of the mean and variance of the lognormal distribution and the Arrhenius law.
Finally, the values of $\tau_1$ and $\tau_2$ are computed such that the probability of failure before $\tau_1$ is $20\%$ and between $\tau_1$ and $\tau_2$ is $40\%$.

For each setting, the bias and Mean Square Error (MSE) for $\gamma_0$, $\gamma_1$ and $\sigma$, based on 1000 simulated samples, are estimated and reported in Tables \ref{t1}-\ref{t4}. The corresponding coverage probabilities and interval lengths of the $90\%,95\%$ and $99\%$ approximate and bootstrap confidence (500 bootstrap replications) intervals are also presented in the tables.

\begin{figure}
\begin{center}
\includegraphics[width=0.49\textwidth]{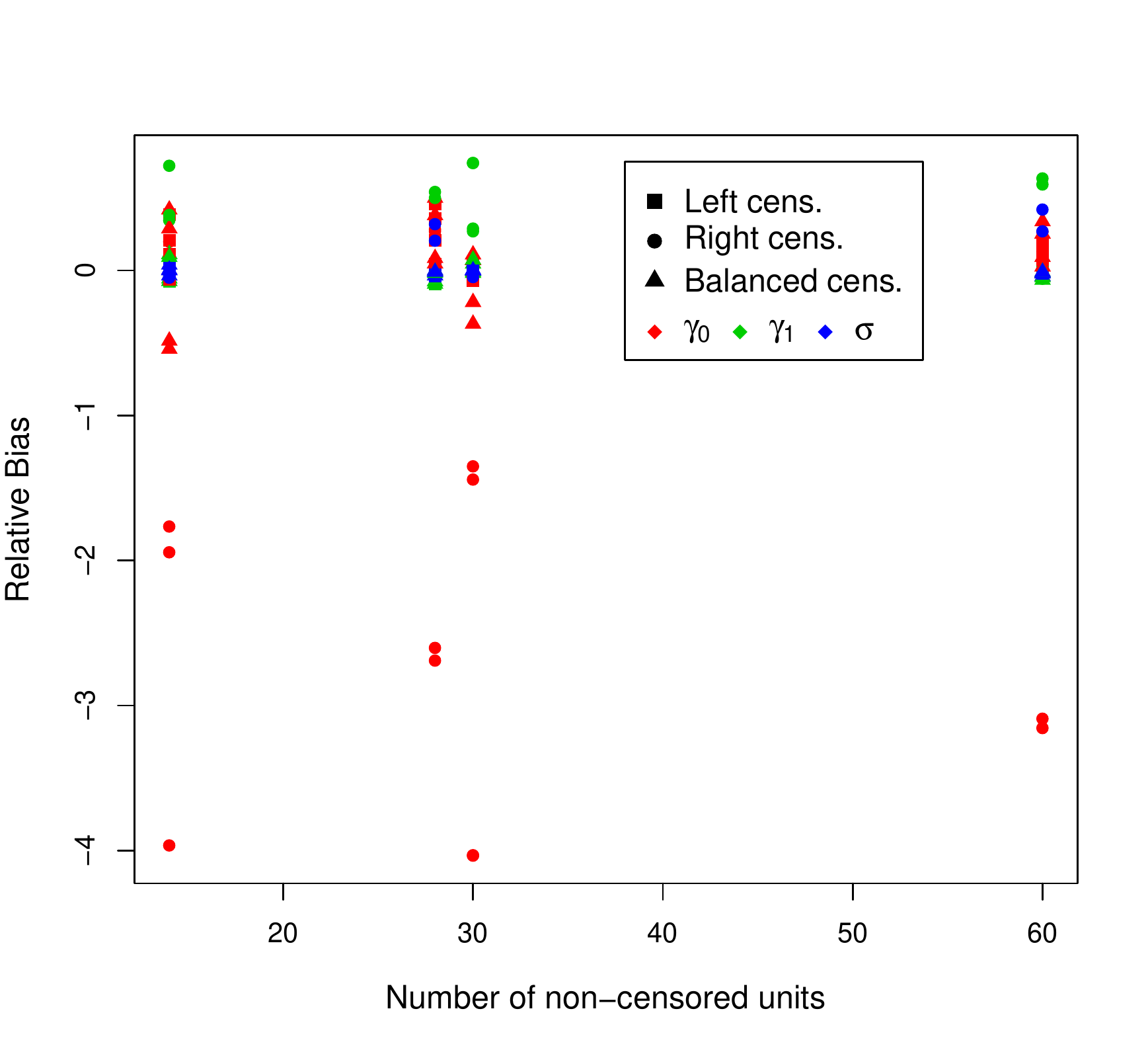}
\includegraphics[width=0.49\textwidth]{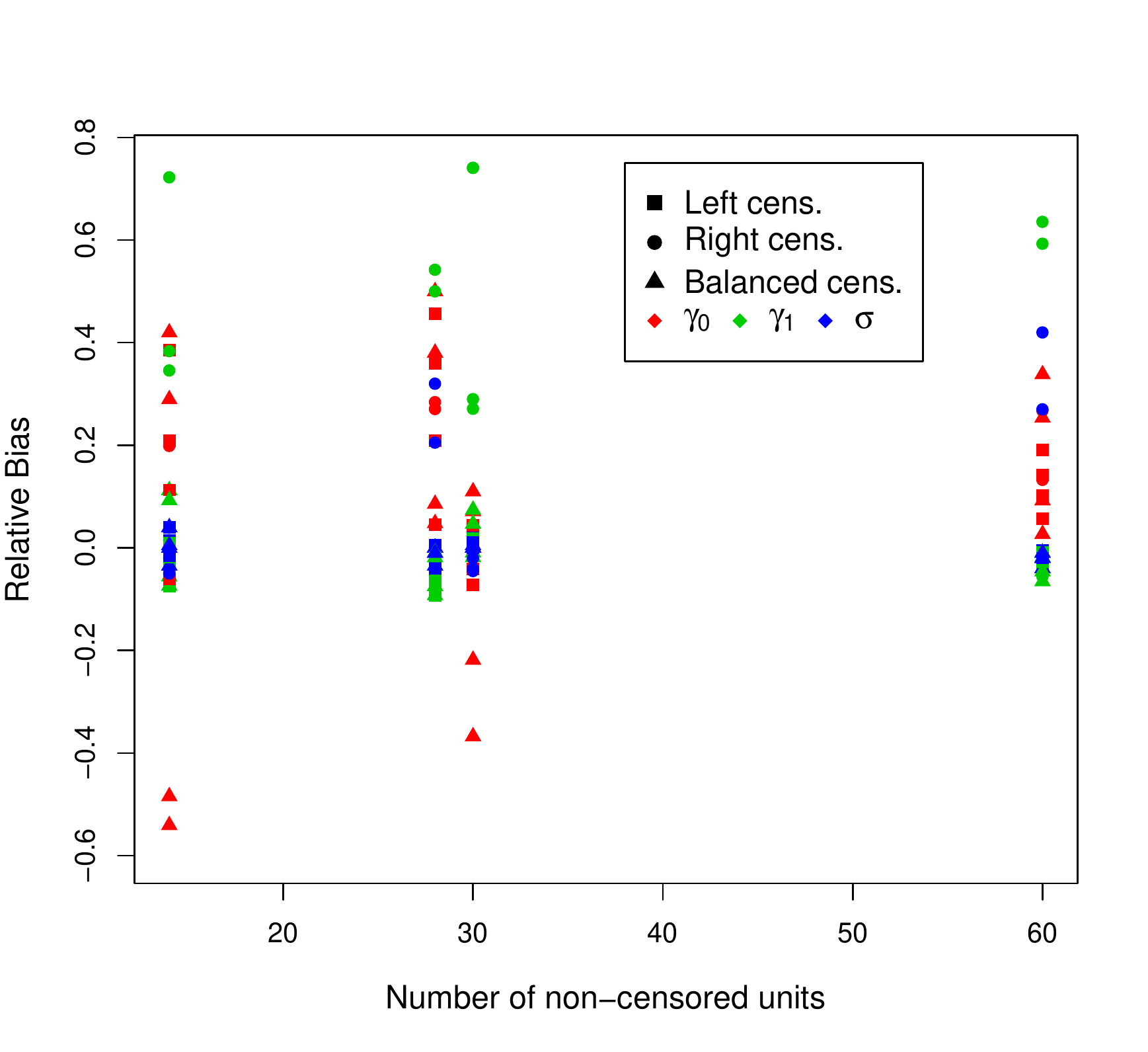}
\caption{\label{fig2}Results of the simulation study: Relative Bias comparison between censoring schemes, number of failures (non-censored units) and parameters. The lower figure shows an enlargement of the upper one for values close to 0.}
\end{center}
\end{figure}

\begin{figure}
\begin{center}
\includegraphics[width=0.49\textwidth]{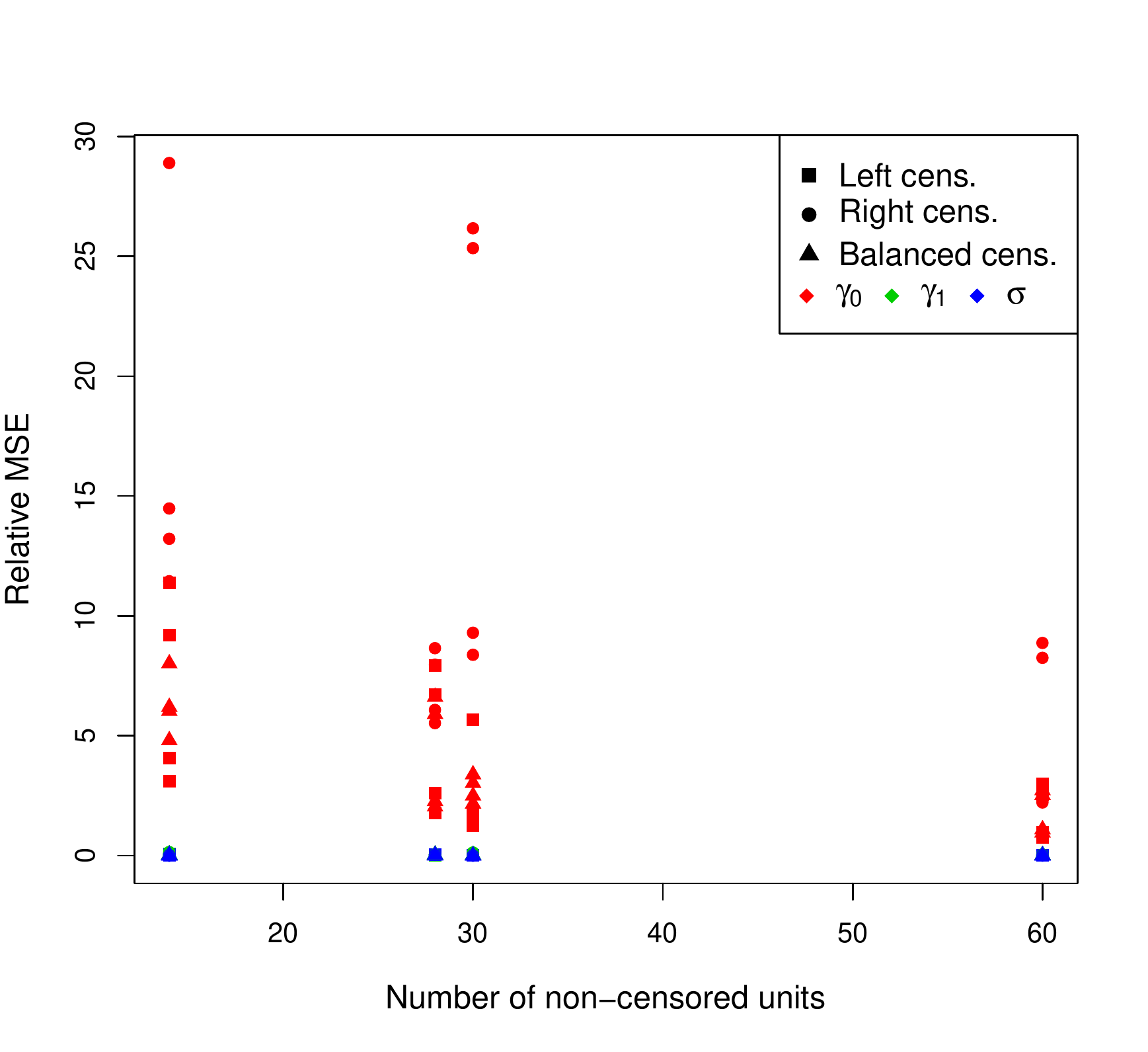}
\includegraphics[width=0.49\textwidth]{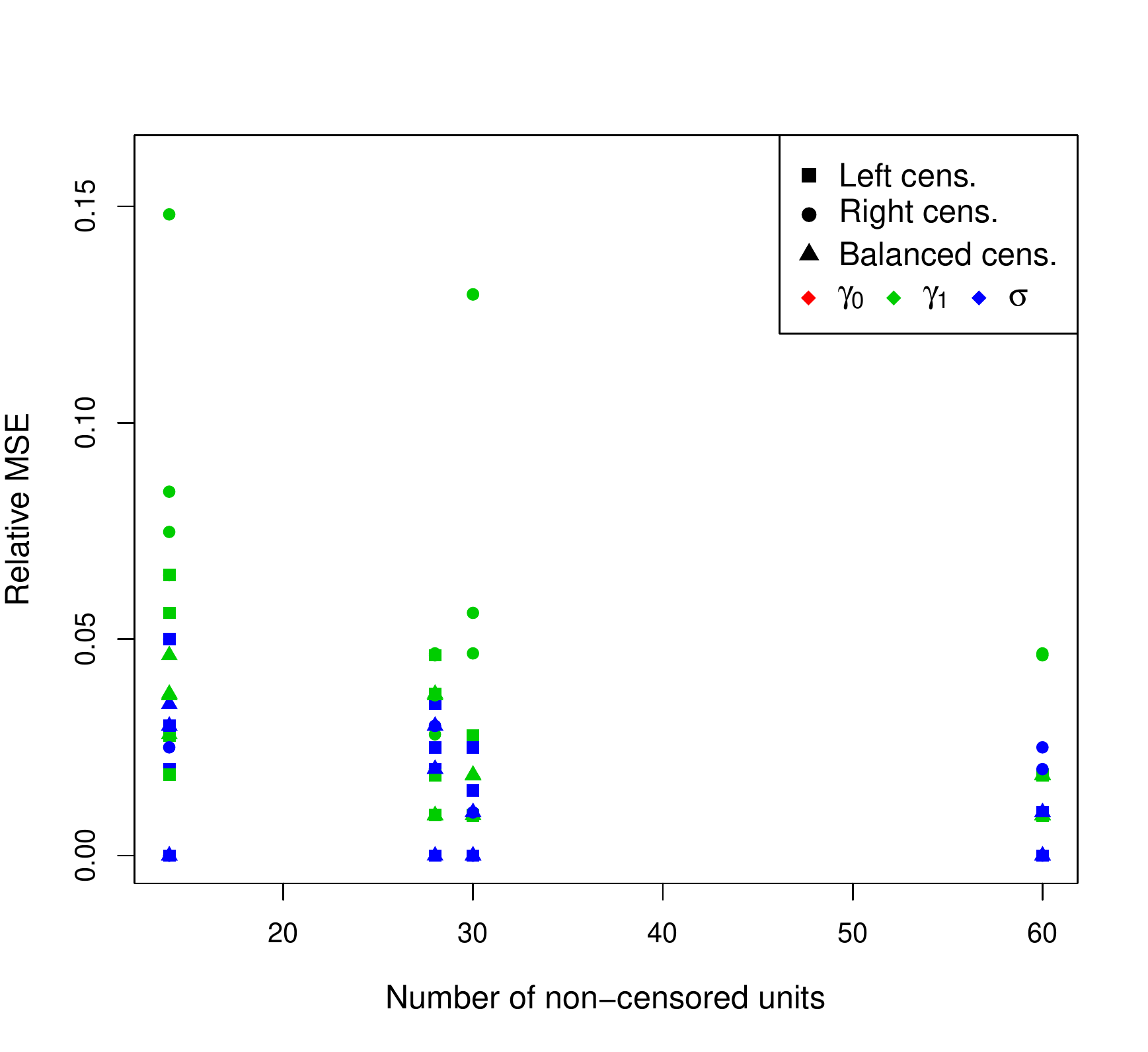}
\caption{\label{fig3}Results of the simulation study: Relative MSE comparison between censoring schemes, number of failures (non-censored units) and parameters. The lower figure shows an enlargement of the upper one for values close to 0.}
\end{center}
\end{figure}

\begin{figure}
\begin{center}
\includegraphics[width=0.49\textwidth]{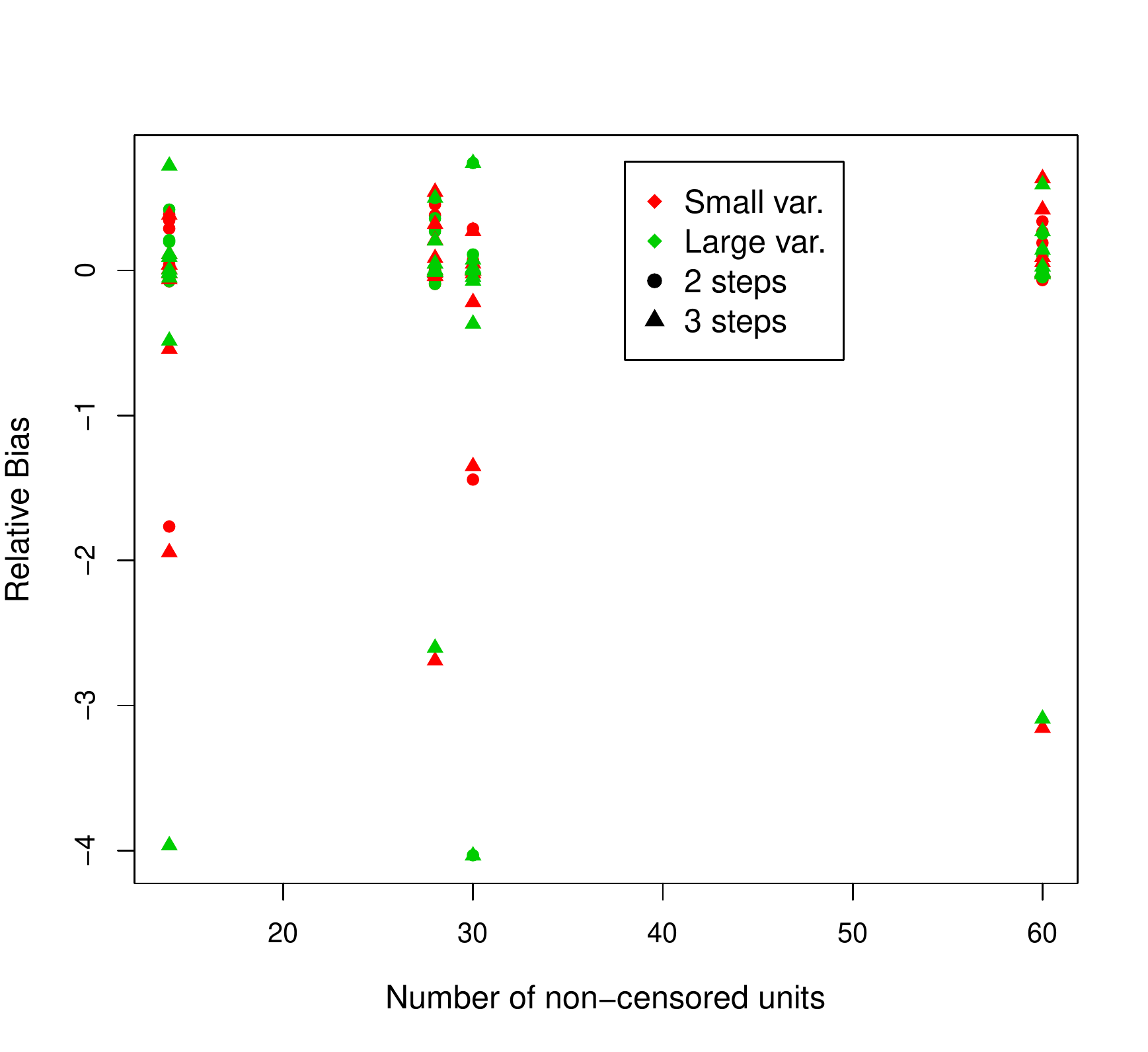}
\includegraphics[width=0.49\textwidth]{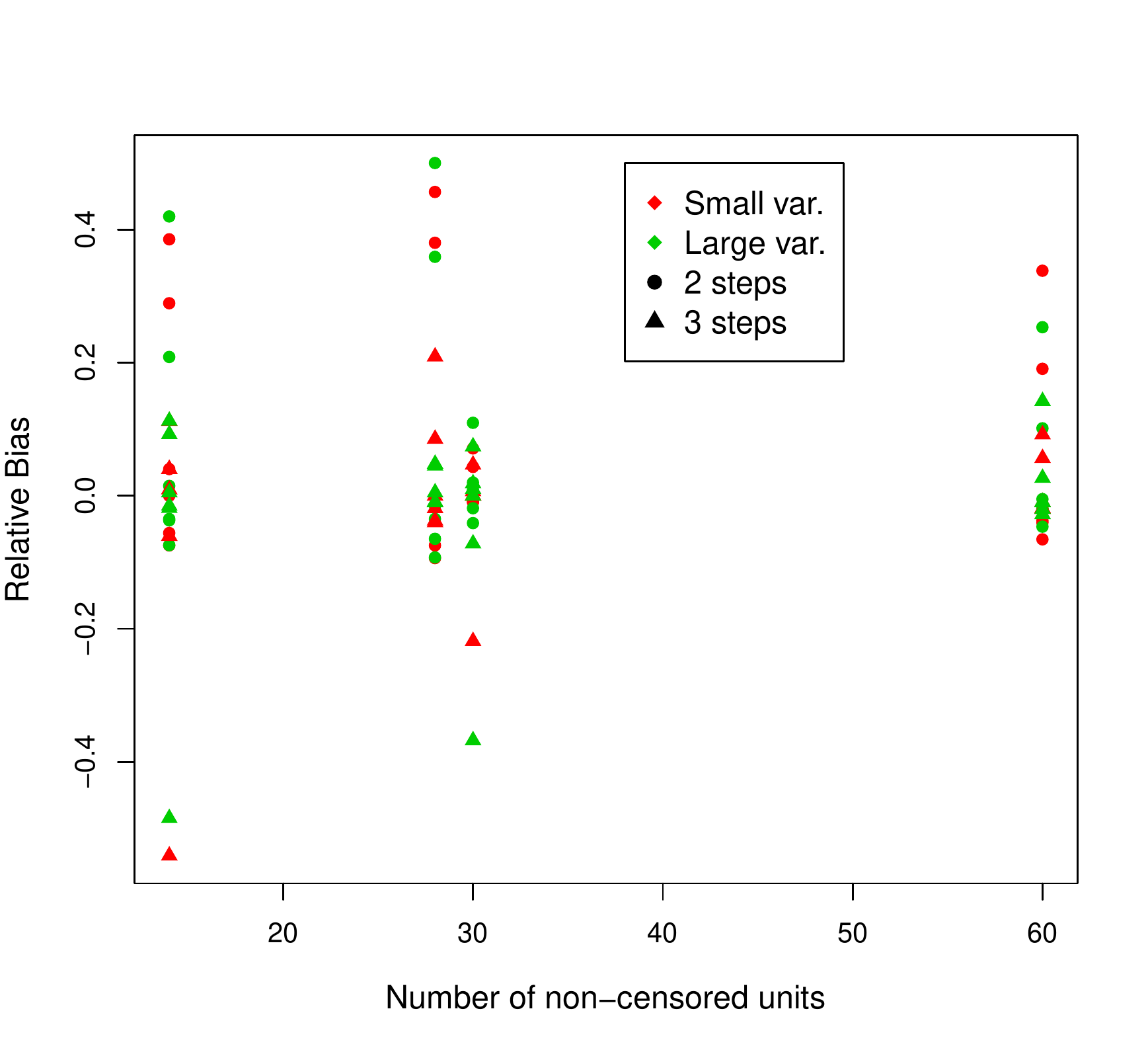}
\caption{\label{fig4}Results of the simulation study: Relative Bias comparison between number of steps, number of failures (non-censored units) and variability level. The lower figure does not contain the classical type-II censoring scheme.}
\end{center}
\end{figure}

\begin{figure}
\begin{center}
\includegraphics[width=0.49\textwidth]{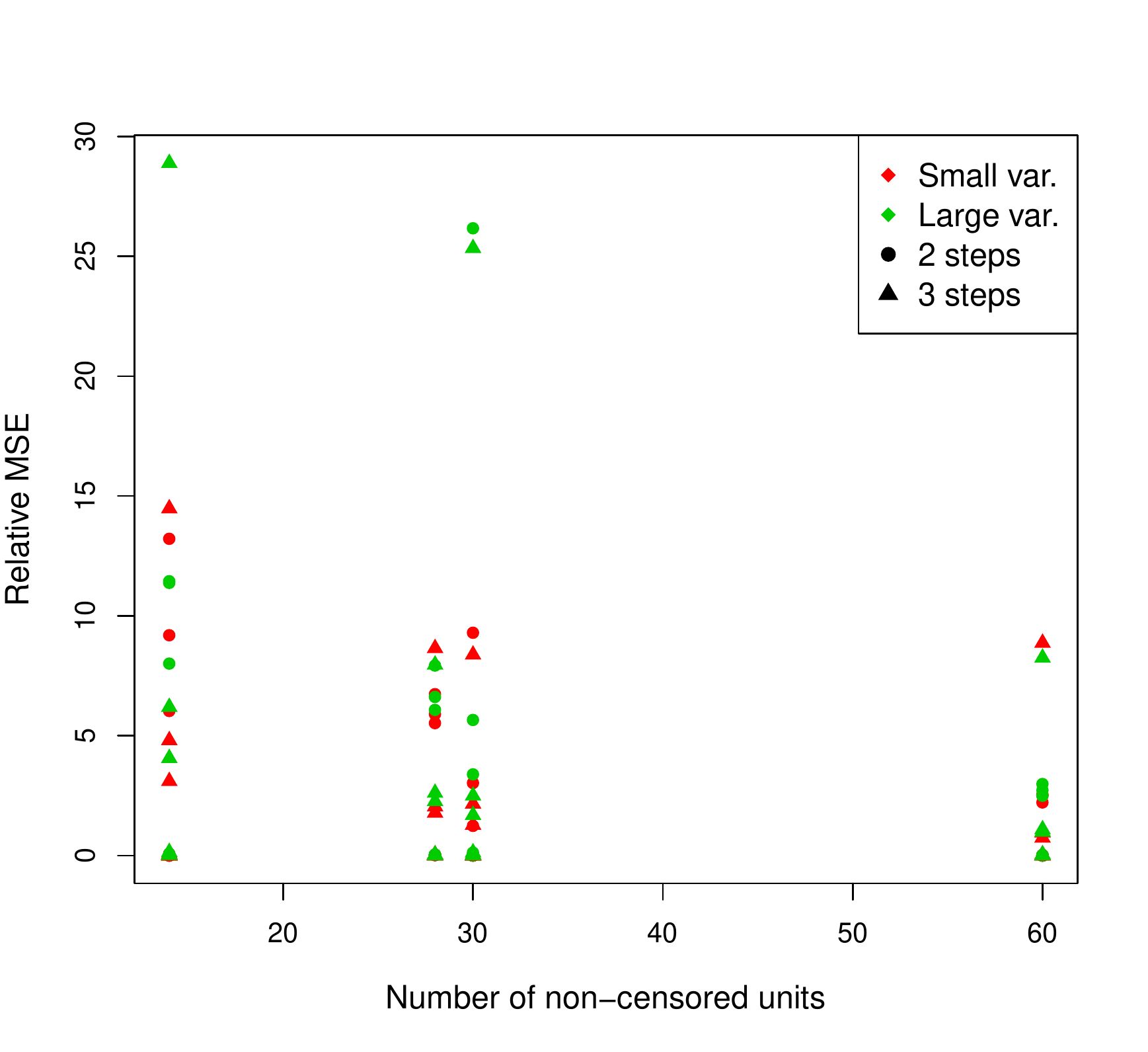}
\includegraphics[width=0.49\textwidth]{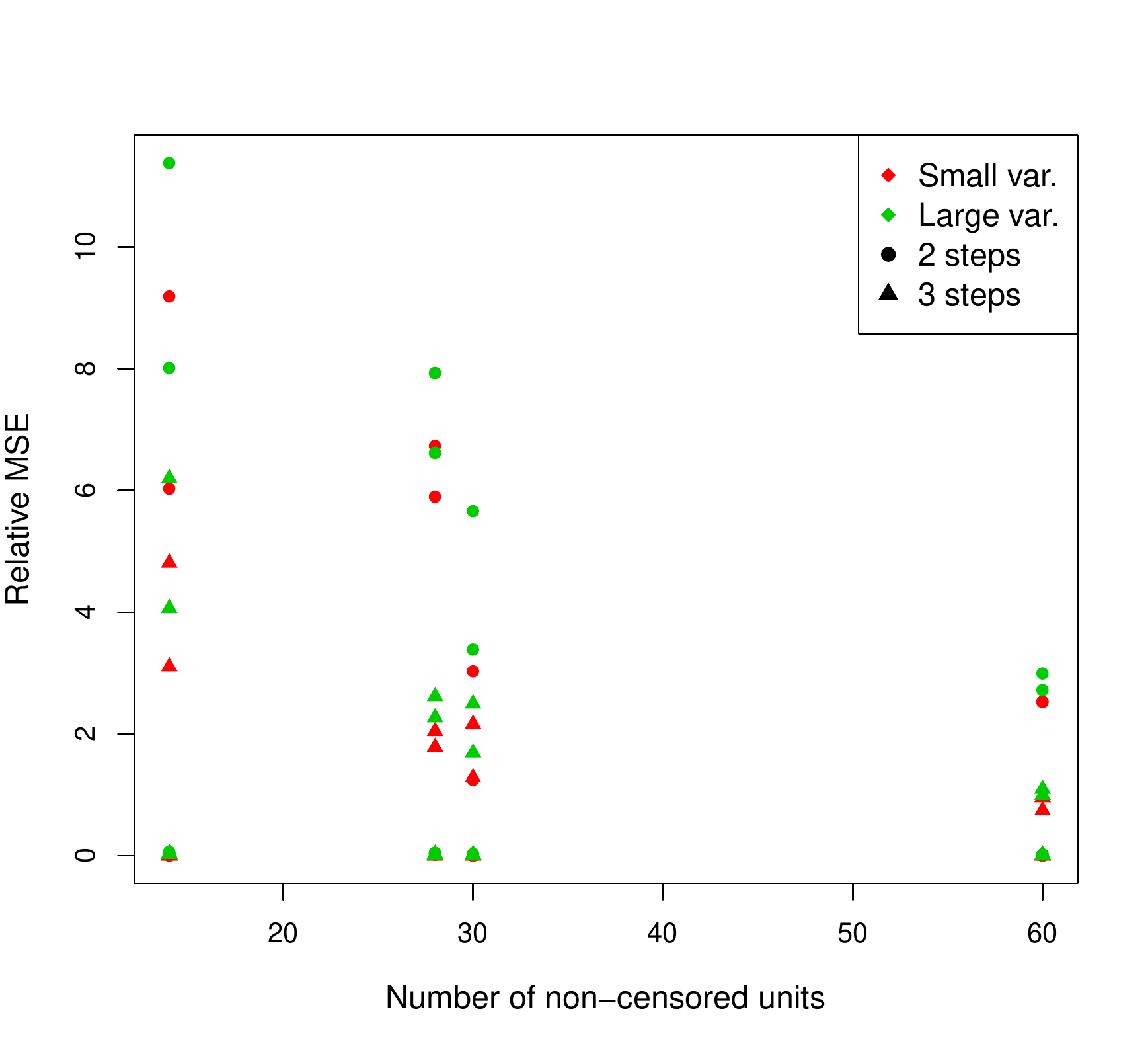}
\caption{\label{fig5}Results of the simulation study: Relative MSE comparison between number of steps, number of failures (non-censored units) and variability level. The lower figure does not contain the classical type-II censoring scheme.}
\end{center}
\end{figure}

Additionally, for the sake of understanding, we have represented the relative bias and relative MSE (which can be obtained by dividing the corresponding bias and MSE by the value of the parameter) in Figures \ref{fig2}-\ref{fig5}, where the x-axis of all the figures contains the number of failures (non-censored units) of each scenario. Note that the first and third ones correspond to high levels of censoring, while the second and fourth ones correspond to low ones. The lower plot of Figures \ref{fig2} and \ref{fig3} shows an enlargement of the upper figures, so as to better evaluate the values close to 0. On the contrary, in the lower plot of Figures \ref{fig4} and \ref{fig5} we have removed the values for the classical Type-II censoring scheme, as it provides much worse results. Through a combination of shapes and colours, Figures \ref{fig2} and \ref{fig3} permit comparing censoring schemes and parameters, whereas Figures \ref{fig4} and \ref{fig5} compare variability level and number of stress levels.

\begin{figure}
\begin{center}
\includegraphics[width=0.49\textwidth]{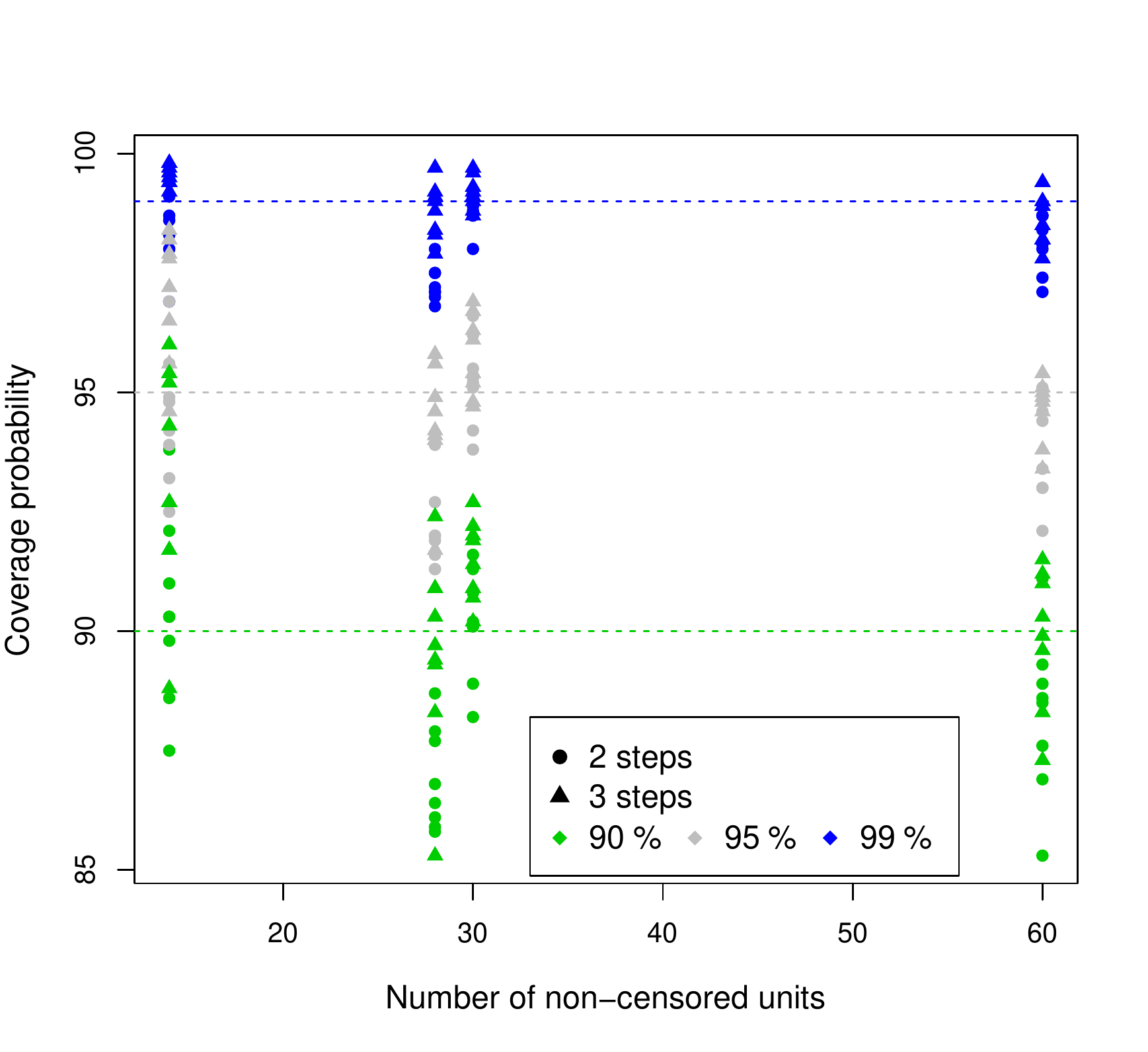}
\includegraphics[width=0.49\textwidth]{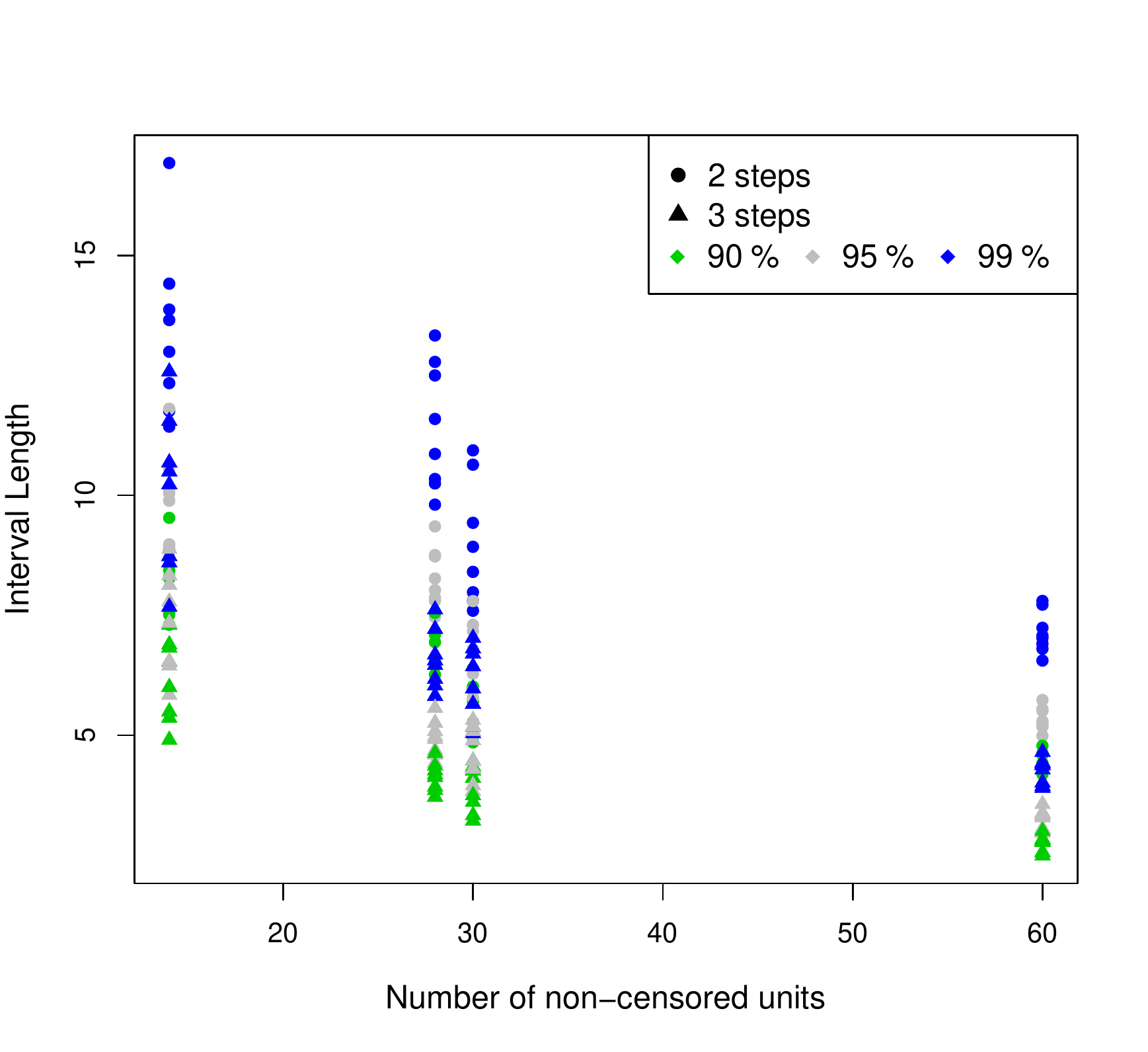}
\caption{\label{fig6}Results of the simulation study: Coverage probabilities (up) and interval lengths (down) comparison between number of steps, nominal confidence level and number of failures (non-censored units). The lower figure does not contain the classical type-II censoring scheme.}
\end{center}
\end{figure}

\begin{figure}
\begin{center}
\includegraphics[width=0.49\textwidth]{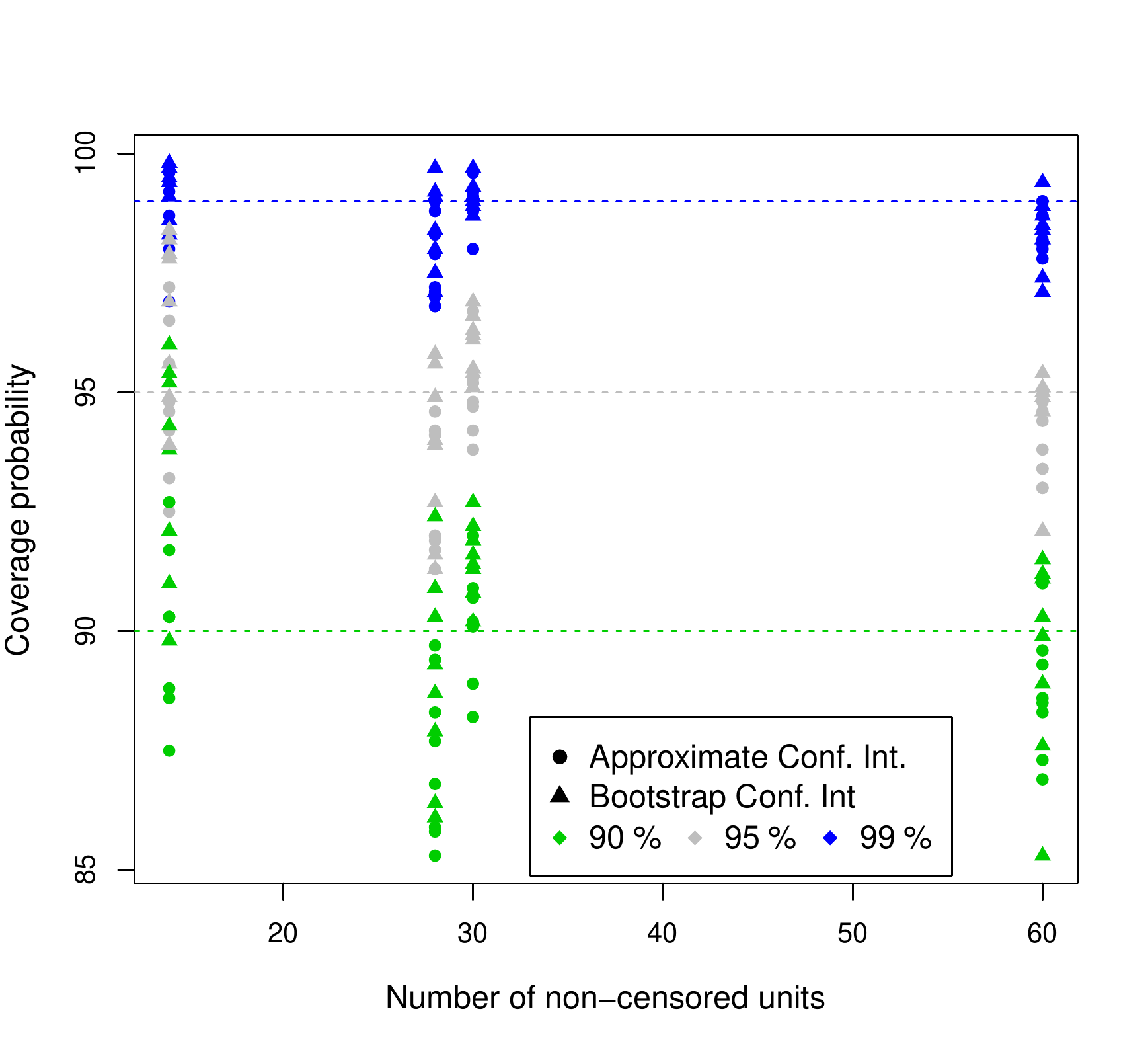}
\includegraphics[width=0.49\textwidth]{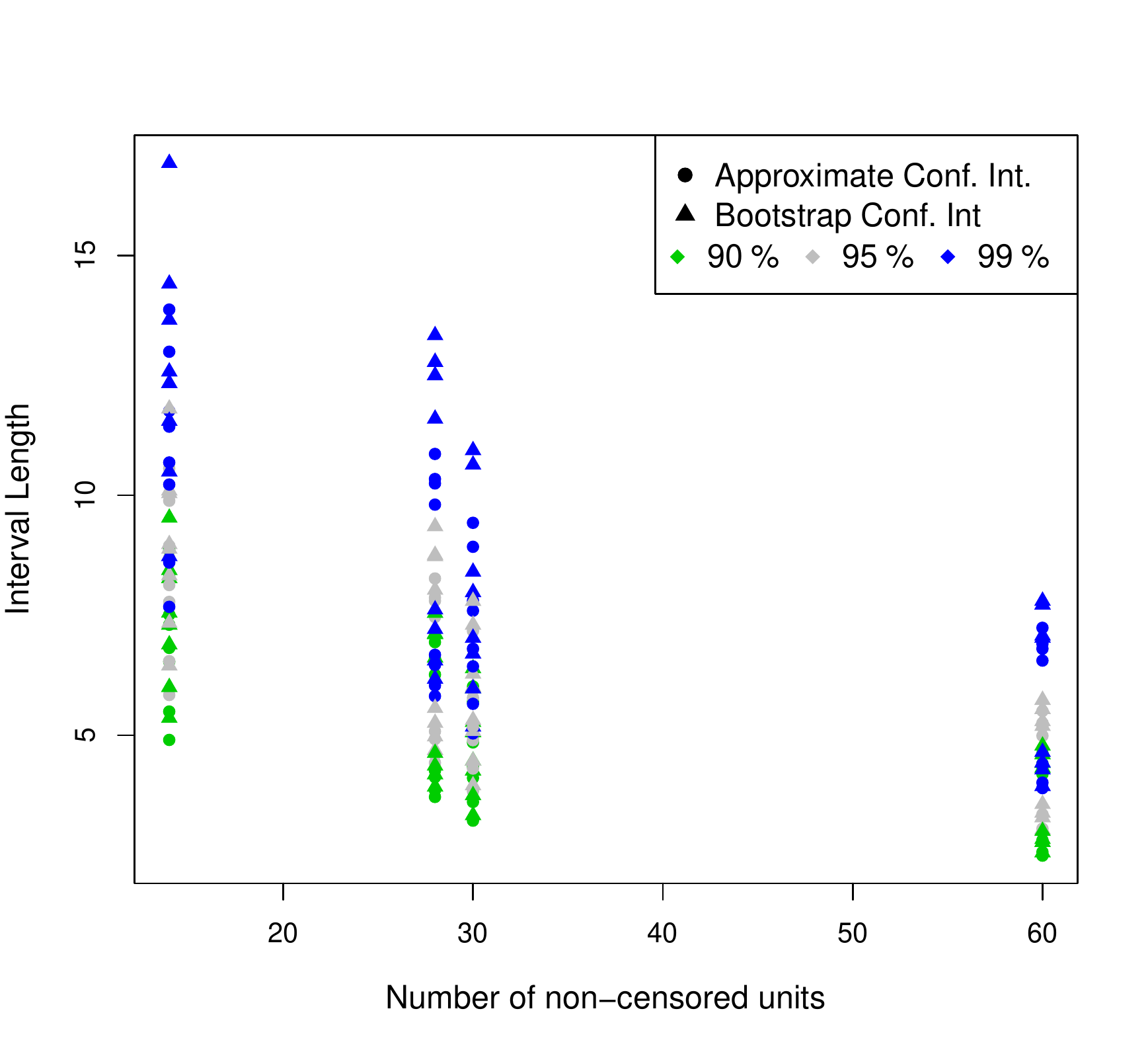}
\caption{\label{fig7}Results of the simulation study: Coverage probabilities (up) and interval lengths (down) comparison between type of interval, nominal confidence level and number of failures (non-censored units). The lower figure does not contain the classical type-II censoring scheme.}
\end{center}
\end{figure}

Regarding confidence intervals, Figures \ref{fig6} and \ref{fig7} contain the coverage probabilities and interval lengths of the CIs associated with parameter $\gamma_0$. Again, the x-axis contains the number of failures (non-censored units) of each scenario. Moreover, we have skipped the complete right censoring schemes because of the bad results it gives (see Tables \ref{t1}-\ref{t4}). Conclusions refering to other parameters are made based on the values of Tables \ref{t1}-\ref{t4}.

For the simulation studies presented here, and others carried out for different values for the parameters (not presented here for the sake of conciseness), we can make the following observations:
\begin{itemize}

\item Classical Type-II censoring (i.e., complete right censoring) leads to worse results in terms of all quantities evaluated in this simulation study, mainly for $k$ equal to $3$. Regarding bias and MSE, Figures \ref{fig2} and \ref{fig3} show that this censoring scheme (represented by a circle) is the one further from 0. For the case of $\gamma_0$, this type of censoring schemes leads to a very large underestimation and small precision of the parameter. $\gamma_0$ is represented in red color in Figures \ref{fig2} and \ref{fig3} . It can be seen that red circles lead to very negative relative bias and very large relative MSE. For that reason, we suggest to avoid this censoring scheme for step-stress models.

  \item Balanced and complete left censoring schemes lead to similar results. Note that in Figures \ref{fig2} and \ref{fig3} neither squares nor triangles outperform the other. In the figures and tables, we have only represented a generic balanced censoring schemes, but there are many other censoring vectors. The ones that lead to better results are those that allocate the censored units among the first failures, as it can be deduced from the results shown in this paper.

  \item When the number of non-censored units gets larger (because of a larger sample size or a smaller proportion of censored units), the estimates are more precise (lower MSE) and therefore, confidence intervals get narrower (smaller interval length). Figures \ref{fig3} and \ref{fig5} shows that the dispersion of the values represented gets smaller when the number of failures gets larger. Moreover, it can be seen that, when the number of non-censored units is similar (second and third values in the x-axis), better results are obtained for smaller proportions of censoring (second value in the x-axis).

  \item For all scenarios, with the exception of the complete right censoring scheme, the estimates for $\gamma_1$ and $\sigma$ are very precise and almost unbiased. Note that, in Figures \ref{fig2} and \ref{fig3}, $\gamma_1$ and $\sigma$ are plotted in green and blue, respectively and those are the predominant colours close to 0 in all plots (specially the lower plot of Figure \ref{fig3}).

  \item An interesting feature of this model is the fact that large variances do not affect the quality of the results: similar levels of accuracy and interval lengths. This fact can be observed in Figures \ref{fig4} and \ref{fig5} where small and large variability are represented in red and green, respectively, and the colours are distributed with no pattern in the plot.

  \item Better results are obtained for 3-steps models than for 2-steps ones, when the remaining parameters are kept constant. In particular, point estimates are more precise, bias are smaller and confidence intervals are narrower. Figures \ref{fig4} and \ref{fig5} represent 2 and 3 steps by circles and triangles, respectively. It can be seen that circles are generally further from $0$ than triangles.

  \item $\gamma_0$ confidence intervals are proportionally wider than those of $\gamma_1$ and $\sigma$ (mainly because of the fact that its corresponding MSE is larger, as seen in Figure \ref{fig3}).

  \item As it can be seen in Tables \ref{t1}-\ref{t4}, all coverage probabilities are close to the nominal ones, except for the complete-right censoring scheme. Moreover, the upper plot of Figure \ref{fig6} shows that coverage probabilities closer to the nominal level are achieved for a larger number of failures. Additionally, we can see that the larger the confidence level, the closer the coverage probabilities.

  \item The number of stress levels (steps) influences interval lengths, but not coverage probabilities, as can be seen in Figure \ref{fig6}, where circles and triangles represent 2 and 3 steps, respectively. Note that the upper plot shows no difference between shapes, but the lower plot has larger values associated with circles (2 steps). The extreme case takes place for a number of failures equal to $60$, where, in Figure \ref{fig6} (lower plot), we can see that the interval lengths grow with confidence level (as expected) but takes larger values for circles than triangles (99\% CIs with 3 steps are as wide as 90\% CIs with 2 steps).

  \item As expected, bootstrap CIs outperform (slightly) approximate CIs in terms of coverage probabilities (see Figure \ref{fig7} and Tables \ref{t1}-\ref{t4}) , specially for small sample sizes. Regarding interval lengths, both methods lead to comparable results for small sample sizes but bootstrap CIs show slightly larger ones for large sample sizes.

  \item Coverage probabilities of the confidence intervals associated with $\sigma$ tend to be smaller than the nominal level.

\end{itemize}

Apart from the previous observations, it can be seen that the proposed inference method behaves well, leading to precise point estimates with small bias and narrow confidence intervals. Better results are obtained when the number of steps is increased and the classical Type-II censoring is avoided.

 \begin{table*}
 \caption{Estimated Bias, MSE, Coverage Probabilities (in percentage) and Interval Lengths (in brackets) of confidence intervals for the three parameters based on $1000$ simulations and $B=500$ replications with $k=2, \gamma_0=0.76, \gamma_1=0.107, \sigma=0.05$ and $\tau_1=95$ for different values of $n$ and $r$ and censoring schemes.}
 \begin{center}
    \resizebox{!}{6.9cm}{
 \renewcommand{\arraystretch}{1.2}
 \renewcommand{\tabcolsep}{1.1mm}
 \begin{tabular}{|c|c|c|c|c|c|c|c|c|c|c|c|}
 \hline
 \multicolumn{6}{|c|}{}&   \multicolumn{2}{c|}{$90\%$ CI} & \multicolumn{2}{c|}{$95\%$ CI} & \multicolumn{2}{c|}{$99\%$ CI}\\
 \hline
 $n$ & $r$ & Cens. scheme & & Bias & MSE & App. & Boot. & App. & Boot. & App. & Boot.\\
  \hline
35 & 14 & $(21,13^\star0)$ & $\gamma_0$ &  0.293 & 6.985 &   87.5 (8.297) &  89.8 (8.440) &   92.5 (9.886) &  93.9 (10.040) &   96.9 (12.993) &  98.3 (13.655)\\
 &  & & $\gamma_1$ &-0.008 & 0.006 &   87.3 (0.235) &  89.7 (0.238) &   92.5 (0.280) &  93.8 (0.284) &   97.0 (0.368) &  98.3 (0.389)\\
&  & & $\sigma$ &0.002 & 0.001 &   80.3 (0.083) &  84.1 (0.077) &   84.0 (0.099) &  88.6 (0.097) &   88.7 (0.130) &  94.5 (0.174)\\
 \cline{3-12}
& & $(13^\star0,21)$ & $\gamma_0$  & -1.342 &10.0411&   78.8 (7.348) &  77.8 (9.608) &   84.9 (8.756) &  80.2 (11.389) &   89.4 (11.507) & 83.0 (15.247)\\
 &  & & $\gamma_1$  &  0.037 &0.008&   79.4 (0.207) &  77.9 (0.269) &   85.1 (0.246) &  80.2 (0.319) &   89.8 (0.324) &  83.0 (0.428)\\
&  & & $\sigma$  & -0.002 &0.000&   80.2 (0.063) &  82.8 (0.066) &   84.7 (0.075) &  88.9 (0.081) &   89.9 (0.099) &  96.1 (0.132)\\
 \cline{3-12}
 & & $(7^\star(1,2))$ & $\gamma_0$  &  0.220 &4.583&   90.3 (7.299) &  93.8 (7.553) &   94.8 (8.697) &  96.9 (8.980) &   98.7 (11.430) &  99.4 (12.337)\\
 &  & & $\gamma_1$  & -0.006 &0.004&   90.4 (0.206) &  93.9 (0.213) &   94.9 (0.245) &  96.9 (0.254) &   98.7 (0.322) &  99.4 (0.351)\\
&  & & $\sigma$  & -0.000 &0.000&   82.4 (0.066) &  88.0 (0.069) &   86.7 (0.079) &  91.6 (0.087) &   91.5 (0.104) &  95.8 (0.167)\\
  \cline{2-12}
& 28 & $(7,27^\star0)$ & $\gamma_0$ &0.347 & 5.115 &   86.8 (6.602) &  88.7 (7.116) &   92.0 (7.867) &  93.9 (8.728) &   96.8 (10.339) &  97.5 (12.779)\\
  &  & & $\gamma_1$ &-0.010 & 0.004 &   86.9 (0.187) &  88.7 (0.201) &   92.0 (0.222) &  93.9 (0.247) &   96.8 (0.292) &  97.5 (0.365)\\
 &  & & $\sigma$ &
0.000 & 0.001 &   79.4 (0.068) &  83.3 (0.071) &   83.6 (0.081) &  89.1 (0.094) &   88.9 (0.106) &  95.5 (0.190)\\
  \cline{3-12}
& & $(27^\star0,7)$ & $\gamma_0$ & 0.216 & 4.204 &   87.6 (6.354) &  90.3 (6.618) &   93.6 (7.571) &  94.7 (7.982) &   98.2 (9.950) &  98.8 (11.358)\\
&  & & $\gamma_1$ &-0.006 & 0.003 &   87.8 (0.179) &  90.1 (0.187) &   93.8 (0.214) &  94.7 (0.226) &   98.2 (0.281) &  98.8 (0.324)\\
&  & & $\sigma$ &-0.000 & 0.001 &   81.7 (0.063) &  84.2 (0.065) &   84.6 (0.076) &  89.6 (0.083) &   90.3 (0.099) &  95.9 (0.166)\\
 \cline{3-12}
& & $(7^\star(0,0,1,0))$ & $\gamma_0$ & 0.289 & 4.484 &   85.8 (6.261) &  86.1 (6.616) &   91.9 (7.460) &  91.3 (8.028) &   97.1 (9.804) &  97.1 (11.594)\\
&  & & $\gamma_1$ &-0.008 & 0.004 &   85.8 (0.177) &  86.2 (0.187) &   92.0 (0.211) &  91.3 (0.227) &   97.0 (0.277) &  97.1 (0.330)\\
&  & & $\sigma$ &0.000 & 0.001 &   79.6 (0.063) &  82.4 (0.065) &   83.9 (0.076) &  88.4 (0.083) &   89.5 (0.099) &  94.8 (0.164)\\
 \hline
75 & 30 & $(45,29^\star0)$ & $\gamma_0$ & 0.033 &   0.947 &  90.1 (5.700) &  90.8 (5.974) &  94.2 (6.792) &  95.5 (7.300) &  98.8 (8.927) &  98.7 (10.637)\\
&  & & $\gamma_1$ &-0.001 &  0.001 &  90.0 (0.162) &  90.7 (0.169) &  94.3 (0.192) &  95.5 (0.207) &  98.8 (0.253) &  98.7 (0.304)\\
 &  & & $\sigma$ &0.000 &  0.000 &  87.4 (0.058) &  90.0 (0.060) &  90.2 (0.069) &  93.6 (0.079) &  94.2 (0.090) &  97.3 (0.148)\\
  \cline{3-12}
& & $(29^\star0,45)$ & $\gamma_0$ &-1.096 &  7.065 &   76.2 (4.831) &  78.5 (7.574) &   81.7 (5.756) &  81.3 (8.787) &   84.2 (7.565) &  83.0 (11.105)\\
&  & & $\gamma_1$ &0.031 &  0.006 &   76.3 (0.136) &  78.5 (0.212) &   81.7 (0.162) &  81.4 (0.245) &   84.2 (0.213) &  83.0 (0.310)\\
 &  & & $\sigma$ &0.000 & 0.000 &   87.2 (0.041) &  87.6 (0.042) &   90.5 (0.048) &  92.0 (0.051) &   94.1 (0.064) &  97.5 (0.072)\\
  \cline{3-12}
  & & $(15^\star(1,2))$ & $\gamma_0$ &0.054 & 2.301 &   90.1 (4.848) &  91.6 (5.060) &   95.2 (5.777) &  96.6 (6.005) &   99.1 (7.592) &  99.1 (7.978)\\
&  & & $\gamma_1$ &-0.001 & 0.002 &   90.4 (0.137) &  91.7 (0.143) &   95.2 (0.163) &  96.9 (0.169) &   99.2 (0.214) &  99.1 (0.225)\\
 &  & & $\sigma$ &0.000 & 0.000 &   88.2 (0.044) &  90.8 (0.046) &   92.3 (0.052) &  95.4 (0.056) &   96.8 (0.069) &  98.6 (0.083)\\
  \cline{2-12}
 & 60 & $(15,59^\star0)$ & $\gamma_0$ & 0.145 & 1.922 &   88.6 (4.410) &  88.9 (4.598) &   93.4 (5.254) &  94.6 (5.542) &   98.7 (6.905) &  98.7 (7.720)\\
&  & & $\gamma_1$ &-0.004 & 0.002 &   88.6 (0.124) &  88.8 (0.130) &   93.4 (0.148) &  94.6 (0.156) &   98.7 (0.195) &  98.7 (0.219)\\
 &  & & $\sigma$ & -0.001 & 0.000 &   85.4 (0.043) &  86.8 (0.045) &   88.9 (0.052) &  91.8 (0.056) &   93.5 (0.068) &  96.8 (0.088)\\
  \cline{3-12}
    & & $(59^\star0,15)$ & $\gamma_0$ &  0.203 & 1.683 &   88.1 (4.145) &  89.2 (4.286) &   93.4 (4.939) &  94.8 (5.134) &   98.6 (6.491) &  98.2 (6.917)\\
&  & & $\gamma_1$ &-0.006 & 0.001 &   88.1 (0.117) &  89.2 (0.121) &   93.5 (0.139) &  94.9 (0.145) &   98.6 (0.183) &  98.2 (0.196)\\
&  & & $\sigma$ &-0.001 & 0.000 &   84.7 (0.040) &  86.0 (0.041) &   88.4 (0.048) &  91.7 (0.050) &   93.8 (0.063) &  96.5 (0.073)\\
     \cline{3-12}
    & & $(15^\star(0,0,1,0))$ & $\gamma_0$ & 0.257 & 1.919 &   86.9 (4.187) &  85.3 (4.318) &   93.0 (4.989) &  91.1 (5.188) &   98.0 (6.557) &  97.4 (7.028)\\
&  & & $\gamma_1$ &-0.007 & 0.002 &   86.7 (0.118) &  85.3 (0.122) &   92.8 (0.141) &  91.1 (0.146) &   98.0 (0.185) &  97.4 (0.199)\\
&  & & $\sigma$ &-0.002 & 0.000 &   82.5 (0.040) &  81.2 (0.041) &   87.4 (0.048) &  88.0 (0.050) &   93.0 (0.063) &  95.1 (0.073)\\
 \hline
 \end{tabular}}
 \end{center}
 \label{t1}
 \end{table*}

 \begin{table*}
 \caption{Estimated Bias, MSE, Coverage Probabilities (in percentage) and Interval Lengths (in brackets) of confidence intervals for the three parameters based on $1000$ simulations and $B=500$ replications with $k=2, \gamma_0=0.71, \gamma_1=0.108, \sigma=0.2$ and $\tau_1=83$ for different values of $n$ and $r$ and censoring schemes.}
 \begin{center}
    \resizebox{!}{6.9cm}{
 \renewcommand{\arraystretch}{1.2}
 \renewcommand{\tabcolsep}{1.1mm}
 \begin{tabular}{|c|c|c|c|c|c|c|c|c|c|c|c|}
 \hline
 \multicolumn{6}{|c|}{} & \multicolumn{2}{c|}{$90\%$ CI} & \multicolumn{2}{c|}{$95\%$ CI} & \multicolumn{2}{c|}{$99\%$ CI}\\
 \hline
 $n$ & $r$ & Cens. scheme & & Bias & MSE & App. & Boot. & App. & Boot. & App. & Boot.\\
  \hline
35 & 14 & $(21,13^\star0)$ & $\gamma_0$ &0.148 & 8.079 &   88.6 (8.859) &  91.0 (9.531) &   93.2 (10.557) &  94.9 (11.805) &   98.0 (13.874) &  99.1 (16.925)\\
&  & & $\gamma_1$ &-0.004 & 0.007 &   88.4 (0.259) &  91.0 (0.278) &   93.2 (0.309) &  95.0 (0.345) &   98.0 (0.406) &  99.0 (0.499)\\
&  & & $\sigma$ &0.003 & 0.010 &   82.2 (0.285) &  86.3 (0.286) &   85.6 (0.339) &  91.5 (0.365) &   90.4 (0.446) &  97.4 (0.58)\\
\cline{3-12}
  & & $(13^\star0,21)$ & $\gamma_0$ &0.141 & 8.129 &   58.5 (7.297) &  50.3 (11.654) &   63.5 (8.695) &  55.0 (14.122) &   73.4 (11.428) &  601 (19.934)\\
&  & & $\gamma_1$ &-0.003 & 0.007 &   60.0 (0.210) &  50.8 (0.327) &   64.5 (0.250) &  55.1 (0.397) &   74.9 (0.329) &  60.5 (0.566)\\
&  & & $\sigma$ &0.003 & 0.010 &   77.8 (0.212) &  78.7 (0.221) &   81.9 (0.252) &  85.1 (0.275) &   88.3 (0.332) &  93.4 (0.458)\\
     \cline{3-12}
& & $(7^\star(1,2))$ & $\gamma_0$ &0.298 & 5.690 &   90.3 (7.511) &  92.1 (8.269) &   94.2 (8.950) &  95.6 (10.097) &   98.3 (11.762) &  98.6 (14.409)\\
&  & & $\gamma_1$ &-0.008 & 0.005 &   90.5 (0.217) &  92.1 (0.239) &   94.5 (0.259) &  95.7 (0.293) &   98.2 (0.340) &  98.7 (0.423)\\
 &  & & $\sigma$ &-0.007 & 0.007 &   81.5 (0.229) &  86.2 (0.245) &   86.0 (0.272) &  90.2 (0.311) &   91.0 (0.358) &  962. (0.515)\\
  \cline{2-12}
 & 28 & $(7,27^\star0)$ & $\gamma_0$ & 0.255 & 5.630 &   85.9 (6.939) &  87.9 (7.548) &   91.3 (8.268) &  92.7 (9.352) &   97.2 (10.866) &  98.0 (13.334)\\
 &  & & $\gamma_1$ &-0.007 & 0.005 &   86.1 (0.202) &  88.2 (0.220) &   91.3 (0.241) &  92.6 (0.273) &   97.2 (0.316) &  98.1 (0.391)\\
 &  & & $\sigma$ &-0.003 & 0.007 &   81.8 (0.230) &  84.5 (0.240) &   85.0 (0.274) &  91.3 (0.300) &   90.9 (0.360) &  97.3 (0.446)\\
  \cline{3-12}
 & & $(27^\star0,7)$ & $\gamma_0$ &  0.192 & 4.318 &   89.5 (6.459) &  91.1 (7.111) &   93.0 (7.696) &  94.2 (8.781) &   98.3 (10.115) &  98.5 (12.976)\\
 &  & & $\gamma_1$ &-0.005 & 0.004 &   89.2 (0.188) &  91.2 (0.207) &   93.1 (0.224) &  94.0 (0.256) &   98.0 (0.294) &  98.4 (0.382)\\
&  & & $\sigma$ &-0.003 & 0.006 &   84.7 (0.224) &  86.9 (0.236) &   87.7 (0.267) &  91.1 (0.298) &   92.2 (0.351) &  97.1 (0.481)\\
\cline{3-12}
 & & $(7^\star(0,0,1,0))$ & $\gamma_0$ &0.355 & 4.699 &   87.7 (6.547) &  86.4 (7.101) &   92.0 (7.801) &  91.6 (8.756) &   97.0 (10.252) &  97.5 (12.500)\\
 &  & & $\gamma_1$ &-0.010 & 0.004 &   87.9 (0.190) &  86.0 (0.206) &   91.8 (0.226) &  91.6 (0.255) &   96.9 (0.297) &  97.4 (0.366)\\
&  & & $\sigma$ &-0.007 & 0.006 &   79.0 (0.216) &  80.9 (0.224) &   85.2 (0.258) &  87.9 (0.281) &   90.0 (0.338) &  95.1 (0.424)\\
 \hline
75 & 30 & $(45,29^\star0)$ & $\gamma_0$ & -0.029 & 4.018 &   88.2 (6.018) &  90.2 (6.391) &   93.8 (7.171) &  95.1 (7.796) &   98.0 (9.424) &  98.9 (10.937)\\
  &  & & $\gamma_1$ &0.001 & 0.003 &   88.5 (0.176) &  89.8 (0.187) &   93.8 (0.209) &  95.2 (0.229) &   98.0 (0.275) &  99.0 (0.323)\\
   &  & & $\sigma$ &0.004 & 0.005 &   85.6 (0.195) &  89.1 (0.207) &   90.0 (0.233) &  93.7 (0.257) &   94.5 (0.306) &  97.5 (0.384)\\
  \cline{3-12}
& & $(29^\star0,45)$ & $\gamma_0$ &  -2.863 & 18.579 &   51.9 (4.810) &  50.3 (8.964) &   55.4 (5.732) &  53.9 (10.427) &   58.3 (7.533) &  57.0 (13.392)\\
&  & & $\gamma_1$ &0.080 & 0.014 &   52.1 (0.138) &  50.4 (0.250) &   55.5 (0.165) &  53.8 (0.291) &   59.0 (0.217) &  57.0 (0.374)\\
 &  & & $\sigma$ &-0.008 & 0.002 &   83.7 (0.142) &  83.5 (0.146) &   88.4 (0.169) &  89.0 (0.177) &   93.2 (0.223) &  95.2 (0.245)\\
\cline{3-12}
 & & $(15^\star(1,2))$ & $\gamma_0$ &0.078 & 2.404 &   88.9 (4.982) &  91.3 (5.272) &   95.3 (5.937) &  96.2 (6.280) &   99.2 (7.802) &  99.1 (8.403)\\
&  & & $\gamma_1$ &-0.002 & 0.002 &   89.4 (0.144) &  92.0 (0.152) &   95.5 (0.171) &  96.1 (0.182) &   99.3 (0.225) &  99.2 (0.244)\\
 &  & & $\sigma$ &0.001 & 0.002 &   90.7 (0.157) &  92.6 (0.166) &   93.7 (0.187) &  96.0 (0.200) &   97.4 (0.246) &  98.6 (0.280)\\
  \cline{2-12}
& 60 & $(15,59^\star0)$ & $\gamma_0$ &  0.072 & 2.124 &   89.3 (4.621) &  91.1 (4.779) &   94.4 (5.506) &  95.1 (5.732) &   98.7 (7.236) &  98.4 (7.798)\\
 &  & & $\gamma_1$ &-0.002 & 0.002 &   89.5 (0.134) &  90.7 (0.139) &   94.4 (0.160) &  94.9 (0.167) &   98.7 (0.210) &  98.5 (0.228)\\
  &  & & $\sigma$ &-0.001 & 0.002 &   87.0 (0.155) &  88.2 (0.159) &   91.8 (0.185) &  93.3 (0.192) &   95.7 (0.243) &  97.8 (0.267)\\
  \cline{3-12}
& & $(59^\star0,15)$ & $\gamma_0$ &0.094 & 1.773 &   88.7 (4.286) &  89.1 (4.404) &   94.6 (5.107) &  95.9 (5.277) &   99.3 (6.712) &  98.6 (7.11)\\
 &  & & $\gamma_1$ &-0.003 & 0.001 &   89.0 (0.124) &  88.9 (0.128) &   94.6 (0.148) &  95.6 (0.153) &   99.3 (0.195) &  98.6 (0.207)\\
&  & & $\sigma$ &-0.002 & 0.002 &   85.2 (0.148) &  86.7 (0.151) &   89.5 (0.176) &  91.7 (0.182) &   95.4 (0.232) &  97.3 (0.254)\\
   \cline{3-12}
 & & $(15^\star(0,0,1,0))$ & $\gamma_0$ &0.180 & 1.931 &   88.5 (4.344) &  87.6 (4.416) &   93.0 (5.177) &  92.1 (5.289) &   98.1 (6.803) &  97.1 (7.074)\\
 &  & & $\gamma_1$ &-0.005 & 0.002 &   88.5 (0.126) &  87.1 (0.128) &   93.0 (0.150) &  92.0 (0.153) &   97.9 (0.197) &  97.0 (0.206)\\
&  & & $\sigma$ &-0.004 & 0.002 &   83.4 (0.144) &  82.3 (0.144) &   88.8 (0.172) &  89.5 (0.174) &   93.8 (0.226) &  95.5 (0.238)\\
 \hline
 \end{tabular}}
 \end{center}
 \label{t2}
 \end{table*}

 \begin{table*}
 \caption{Estimated Bias, MSE, Coverage Probabilities (in percentage) and Interval Lengths (in brackets) of confidence intervals for the three parameters based on $1000$ simulations and $B=500$ replications with $k=3, \gamma_0=0.76, \gamma_1=0.107, \sigma=0.05, \tau_1=95$ and $\tau_2=97.5$ for different values of $n$ and $r$ and censoring schemes.}
\begin{center}
    \resizebox{!}{6.9cm}{
 \renewcommand{\arraystretch}{1.2}
 \renewcommand{\tabcolsep}{1.1mm}
 \begin{tabular}{|c|c|c|c|c|c|c|c|c|c|c|c|}
 \hline
 \multicolumn{6}{|c|}{}&  \multicolumn{2}{c|}{$90\%$ CI} & \multicolumn{2}{c|}{$95\%$ CI} & \multicolumn{2}{c|}{$99\%$ CI}\\
 \hline
 $n$ & $r$ & Cens. scheme & & Bias & MSE & App. & Boot. & App. & Boot. & App. & Boot.\\
  \hline
35 & 14 & $(21,13^\star0)$ & $\gamma_0$ &  -0.046  & 2.360 &   88.8 (4.901) &  95.2 (5.359) &   94.6 (5.839) &  97.9 (6.449) &   99.2 (7.674) &  99.8 (8.729)\\
&  & & $\gamma_1$ &0.001  & 0.002 &   89.0 (0.139) &  95.2 (0.152) &   94.9 (0.165) &  97.9 (0.183) &   99.2 (0.217) &  99.8 (0.249)\\
&  & & $\sigma$ &0.002  & 0.000 &   84.7 (0.064) &  89.9 (0.068) &   89.3 (0.076) &  95.0 (0.084) &   94.6 (0.100) &  98.3 (0.132)\\
  \cline{3-12}
& & $(13^\star0,21)$ & $\gamma_0$ &  -1.477& 11.004 &  78.0 (7.302) &  78.0 (9.802) &   83.6 (8.701) &  80.5 (11.697) &   90.1 (11.434) &  82.8 (16.021)\\
&  & & $\gamma_1$ &0.041&  0.009 &   78.2 (0.206) &  78.0 (0.275) &   84.0 (0.245) &  80.5 (0.328) &   90.3 (0.322) &  82.8 (0.454)\\
&  & & $\sigma$ &-0.001&  0.000 &   79.9 (0.062) &  82.6 (0.071) &   83.1 (0.074) &  88.4 (0.092) &   89.0 (0.097) &  95.7 (0.182)\\
\cline{3-12}
& & $(7^\star(1,2))$ & $\gamma_0$ &-0.411 & 3.655 &  92.7 (6.530) &  96.0 (6.890) &   97.2 (7.781) &  98.4 (8.330) &   99.6 (10.226) &  99.8 (11.550)\\
&  & & $\gamma_1$ &0.012 & 0.003 &   92.8 (0.184) &  96.1 (0.195) &   97.2 (0.219) &  98.4 (0.236) &   99.6 (0.288) &  99.8 (0.330)\\
&  & & $\sigma$ &0.002 & 0.000 &   85.1 (0.065) &  90.3 (0.074) &   87.2 (0.077) &  93.9 (0.097) &   92.6 (0.102) &  97.7 (0.178)\\
  \cline{2-12}
 & 28 & $(7,27^\star0)$ & $\gamma_0$ &  0.159 &  1.356 &   88.3 (3.713) &  90.9 (3.925) &   94.2 (4.425) &  95.8 (4.686) &   98.8 (5.815) &  99.7 (6.173)\\
 &  & & $\gamma_1$ &-0.004 &  0.001 &   88.2 (0.105) &  90.8 (0.111) &   94.1 (0.125) &  96.0 (0.133) &   98.8 (0.164) &  99.7 (0.175)\\
  &  & & $\sigma$ &-0.002 &  0.000 &   86.2 (0.052) &  87.5 (0.053) &   89.1 (0.062) &  91.7 (0.063) &   93.7 (0.082) &  96.7 (0.087)\\
  \cline{3-12}
& & $(27^\star0,7)$ & $\gamma_0$ &-2.045 & 6.574 &   61.2 (5.093) &  32.6 (6.827) &   75.0 (6.069) &  47.0 (8.366) &   92.2 (7.976) &  71.4 (11.881)\\
 &  & & $\gamma_1$ &0.058 & 0.005 &   61.4 (0.145) &  32.7 (0.197) &   75.8 (0.173) &  47.3 (0.242) &   93.0 (0.227) &  71.4 (0.347)\\
  &  & & $\sigma$ &0.016 & 0.001 &   95.8 (0.084) &  66.6 (0.146) &   98.0 (0.100) &  81.1 (0.186) &   98.9 (0.131) &  95.4 (0.305)\\
  \cline{3-12}
& & $(7^\star(0,0,1,0))$ & $\gamma_0$ &0.065 &  1.554 &   85.3 (3.854) &  89.3 (4.180) &   91.7 (4.592) &  94.0 (4.971) &   97.9 (6.035) &  99.2 (6.555)\\
 &  & & $\gamma_1$ &-0.002 &  0.001 &   85.8 (0.109) &  89.3 (0.118) &   91.9 (0.130) &  93.9 (0.140) &   97.9 (0.170) &  99.2 (0.185)\\
  &  & & $\sigma$ &-0.000 &  0.000 &   84.6 (0.051) &  86.0 (0.053) &   87.3 (0.061) &  90.5 (0.063) &   90.8 (0.080) &  95.9 (0.084)\\
 \hline
75 & 30 & $(45,29^\star0)$ & $\gamma_0$ & 0.005 &  0.975 &   90.7 (3.216) &  92.2 (3.329) &   94.7 (3.832) &  96.1 (3.952) &   98.8 (5.036) &  98.7 (5.174)\\
 &  & & $\gamma_1$ &0.000 & 0.001 &   90.7 (0.091) &  92.1 (0.094) &   94.8 (0.108) &  96.1 (0.112) &   98.8 (0.143) &  98.7 (0.147)\\
   &  & & $\sigma$ &0.000 & 0.000 &   88.0 (0.043) &  89.0 (0.044) &   91.7 (0.051) &  93.7 (0.052) &   96.0 (0.067) &  97.1 (0.070)\\
  \cline{3-12}
 & & $(29^\star0,45)$ & $\gamma_0$ & -1.026 & 6.368 &   77.9 (4.849) &  79.2 (7.572) &   81.0 (5.778) &  81.0 (8.776) &   83.3 (7.593) &  83.2 (11.157)\\
&  & & $\gamma_1$ &0.029 & 0.005 &   77.9 (0.136) &  79.2 (0.212) &   81.1 (0.163) &  81.0 (0.245) &   83.6 (0.214) &  83.2 (0.312)\\
&  & & $\sigma$ &-0.001 & 0.000 &   86.6 (0.040) &  87.7 (0.042) &   90.7 (0.048) &  92.0 (0.051) &   94.7 (0.063) &  97.3 (0.075)\\
\cline{3-12}
& & $(15^\star(1,2))$ & $\gamma_0$ &-0.166 & 1.643 &   90.2 (4.107) &  91.4 (4.260) &   94.8 (4.894) &  95.4 (5.083) &   99.2 (6.432) &  99.0 (6.699)\\
&  & & $\gamma_1$ &0.005 & 0.001 &   90.3 (0.116) &  91.5 (0.120) &   94.9 (0.138) &  95.4 (0.143) &   99.2 (0.181) &  99.0 (0.189)\\
&  & & $\sigma$ &0.000 & 0.000 &   88.3 (0.041) &  89.3 (0.043) &   91.6 (0.049) &  93.3 (0.052) &   95.8 (0.065) &  97.7 (0.071)\\
  \cline{2-12}
& 60 & $(15,59^\star0)$ & $\gamma_0$ &  0.043 &  0.564 &   91.0 (2.489) &  91.5 (2.553) &   94.6 (2.966) &  95.4 (3.031) &   98.2 (3.898) &  98.5 (3.937)\\
&  & & $\gamma_1$ &-0.001 &  0.001 &   91.0 (0.070) &  91.3 (0.072) &   94.6 (0.084) &  95.5 (0.086) &   98.3 (0.110) &  98.6 (0.111)\\
&  & & $\sigma$ &-0.001 &  0.000 &   90.6 (0.036) &  90.7 (0.036) &   93.7 (0.042) &  94.3 (0.043) &   97.1 (0.056) &  97.7 (0.056)\\
  \cline{3-12}
 & & $(59^\star0,15)$ & $\gamma_0$ &-2.398 & 6.743 &   23.5 (3.421) &  5.3 (4.608) &   34.7 (4.077) &  8.1 (5.526) &   60.4 (5.357) &  20.6 (7.411)\\
&  & & $\gamma_1$ &0.068 & 0.005 &   23.7 (0.097) &  5.5 (0.133) &   35.3 (0.116) & 8.2 (0.160) &   61.8 (0.152) &  20.6 (0.215)\\
&  & & $\sigma$ &0.021 & 0.001 &   82.1 (0.059) &  24.9 (0.104) &   93.8 (0.070) &  33.9 (0.126) &   99.7 (0.092) &  57.2 (0.174)\\
  \cline{3-12}
& & $(15^\star(0,0,1,0))$ & $\gamma_0$ &0.070 & 0.727 &   88.3 (2.561) &  90.3 (2.770) &   93.4 (3.051) &  95.0 (3.285) &   99.0 (4.010) &  99.4 (4.282)\\
&  & & $\gamma_1$ &-0.002 & 0.001 &   88.6 (0.072) &  90.3 (0.078) &   93.4 (0.086) &  95.0 (0.093) &   99.1 (0.113) &  99.4 (0.121)\\
&  & & $\sigma$ &-0.001 & 0.000 &   88.9 (0.035) &  90.2 (0.035) &   92.2 (0.041) &  93.1 (0.042) &   96.3 (0.054) &  97.2 (0.055)\\
 \hline
 \end{tabular}}
 \end{center}
 \label{t3}
 \end{table*}

  \begin{table*}
 \caption{Estimated Bias, MSE, Coverage Probabilities (in percentage) and Interval Lengths (in brackets) of confidence intervals for the three parameters based on $1000$ simulations and $B=500$ replications with $k=3, \gamma_0=0.71, \gamma_1=0.108, \sigma=0.2, \tau_1=83$ and $\tau_2=92$ for different values of $n$ and $r$ and censoring schemes.}
\begin{center}
    \resizebox{!}{6.9cm}{
 \renewcommand{\arraystretch}{1.2}
 \renewcommand{\tabcolsep}{1.1mm}
 \begin{tabular}{|c|c|c|c|c|c|c|c|c|c|c|c|}
 \hline
 \multicolumn{6}{|c|}{} & \multicolumn{2}{c|}{$90\%$ CI} & \multicolumn{2}{c|}{$95\%$ CI} & \multicolumn{2}{c|}{$99\%$ CI}\\
 \hline
 $n$ & $r$ & Cens. scheme & & Bias & MSE & App. & Boot. & App. & Boot. & App. & Boot.\\
  \hline
35 & 14 & $(21,13^\star0)$ & $\gamma_0$ &  0.080  & 2.887 &   91.7 (5.491) &  94.3 (6.000) &   95.6 (6.543) &  97.8 (7.338) &   99.4 (8.599) &  99.5 (10.494)\\
 &  & & $\gamma_1$ &-0.002  & 0.003 &   92.0 (0.162) &  94.4 (0.177) &   96.1 (0.193) &  97.8 (0.217) &   99.3 (0.254) &  99.4 (0.315)\\
   &  & & $\sigma$ &-0.003  & 0.006 &   85.5 (0.250) &  88.9 (0.258) &   88.3 (0.297) &  93.9 (0.326) &   93.5 (0.391) &  97.8 (0.536)\\
 \cline{3-12}
    & & $(13^\star0,21)$ & $\gamma_0$ &-2.815 & 20.508 &   61.6 (7.365) &  56.7 (11.502) &   67.1 (8.776) &  60.0 (13.977) &   76.2 (11.533) &  64.2 (19.787)\\
 &  & & $\gamma_1$ &0.078 &  0.016 &   62.7 (0.212) &  56.7 (0.323) &   67.6 (0.253) &  60.1 (0.394) &   77.6 (0.332) &  64.8 (0.565)\\
   &  & & $\sigma$ &-0.010 &  0.005 &   81.2 (0.215) &  82.1 (0.229) &   84.1 (0.256) &  87.5 (0.290) &   89.9 (0.337) &  94.3 (0.495)\\
     \cline{3-12}
& & $(7^\star(1,2))$ & $\gamma_0$ &-0.344 & 4.400 &   92.7 (6.822) &  95.4 (7.300) &   96.5 (8.129) &  98.2 (8.882) &   99.4 (10.683) &  99.7 (12.578)\\
&  & & $\gamma_1$ &0.010 & 0.004 &   93.2 (0.198) &  95.4 (0.212) &   96.8 (0.236) &  98.4 (0.259) &   99.4 (0.310) &  99.7 (0.372)\\
&  & & $\sigma$ & 0.001 & 0.006  &   83.5 (0.238) &  87.7 (0.254) &   87.8 (0.283) &  91.7 (0.316) &   92.0 (0.372) &  97.2 (0.506)\\
  \cline{2-12}
  & 28 & $(7,27^\star0)$ & $\gamma_0$ & 0.032 &  1.858 &   89.4 (4.126) &  90.3 (4.364) &   94.6 (4.916) &  94.9 (5.256) &   98.3 (6.461) &  98.4 (7.211)\\
  &  & & $\gamma_1$ &-0.001 &  0.002 &   89.3 (0.121) &  90.5 (0.128) &   94.5 (0.144) &  95.4 (0.154) &   98.3 (0.189) &  98.3 (0.213)\\
&  & & $\sigma$ &0.001 &  0.005 &   85.2 (0.198) &  86.2 (0.201) &   88.8 (0.236) &  91.7 (0.245) &   93.8 (0.310) &  96.2 (0.358)\\
\cline{3-12}
& & $(27^\star0,7)$ & $\gamma_0$ &-1.848 & 5.655 &   69.4 (5.408) &  44.1 (7.117) &   82.0 (6.444) &  59.1 (8.668) &   96.7 (8.469) &  83.6 (12.102)\\
 &  & & $\gamma_1$ &0.054 & 0.005 &   71.3 (0.160) &  44.6 (0.216) &   84.0 (0.191) &  59.8 (0.264) &   97.8 (0.250) &  84.3 (0.370)\\
&  & & $\sigma$ &0.041 & 0.006 &   95.9 (0.270) &  83.6 (0.381) &   97.8 (0.322) &  92.6 (0.467) &   99.1 (0.423) &  98.4 (0.670)\\
\cline{3-12}
& & $(7^\star(0,0,1,0))$ &$\gamma_0$ &0.034 & 1.611 &   89.7 (4.264) &  92.4 (4.627) &   94.1 (5.081) &  95.6 (5.564) &   99.0 (6.678) &  99.1 (7.617)\\
&  & & $\gamma_1$ &-0.001 & 0.001 &   89.7 (0.125) &  92.0 (0.135) &   94.2 (0.148) &  95.7 (0.163) &   98.9 (0.195) &  99.1 (0.224)\\
  &  & & $\sigma$ &-0.002 & 0.004 &   85.7 (0.196) &  88.5 (0.201) &   89.9 (0.233) &  91.9 (0.243) &   94.1 (0.307) &  97.1 (0.352)\\
  \hline
75 & 30 & $(45,29^\star0)$ & $\gamma_0$ &  -0.051 &  1.200 &   90.9 (3.609) &  91.9 (3.746) &   95.2 (4.301) &  96.3 (4.465) &   99.1 (5.652) &  99.3 (5.971)\\
 &  & & $\gamma_1$ &0.002 &  0.001 &   90.9 (0.106) &  92.1 (0.110) &   95.3 (0.126) &  96.5 (0.131) &   99.3 (0.166) &  99.3 (0.177)\\
&  & & $\sigma$ & 0.002 &  0.003 &   88.6 (0.165) &  90.5 (0.169) &   92.3 (0.196) &  95.2 (0.203) &   96.7 (0.258) &  98.4 (0.285)\\
\cline{3-12}
   & & $(29^\star0,45)$ & $\gamma_0$ &-2.865 & 17.989 &   52.2 (4.815) &  48.5 (8.954) &   55.3 (5.737) &  52.7 (10.393) &   57.8 (7.540) &  56.4 (13.358)\\
 &  & & $\gamma_1$ &0.080 &  0.014 &   52.3 (0.138) &  48.7 (0.250) &   55.4 (0.165) &  52.6 (0.290) &   58.0 (0.217) &  56.5 (0.373)\\
&  & & $\sigma$ & -0.009 & 0.002&   84.7 (0.142) &  84.9 (0.147) &   89.5 (0.17) &  89.6 (0.177) &   94.1 (0.223) &  95.7 (0.244)\\
\cline{3-12}
& & $(15^\star(1,2))$ & $\gamma_0$ & -0.261 &  1.774 &   92.0 (4.345) &  92.7 (4.466) &   96.7 (5.177) &  96.9 (5.316) &   99.6 (6.804) &  99.7 (7.027)\\
 &  & & $\gamma_1$ &0.008 &  0.002 &   92.1 (0.126) &  93.1 (0.129) &   96.9 (0.150) &  97.3 (0.154) &   99.6 (0.197) &  99.7 (0.204)\\
&  & & $\sigma$ &-0.000 &  0.002 &   88.4 (0.150) &  90.4 (0.156) &   91.7 (0.179) &  94.6 (0.186) &   96.4 (0.235) &  98.3 (0.251)\\
  \cline{2-12}
 & 60 & $(15,59^\star0)$ & $\gamma_0$ &0.101 &  0.705 &   89.6 (2.791) &  91.2 (2.832) &   94.8 (3.325) &  94.9 (3.380) &   99.0 (4.370) &  98.2 (4.429)\\
  &  & & $\gamma_1$ &-0.003 &  0.001 &   89.6 (0.082) &  90.8 (0.083) &   94.7 (0.097) &  94.9 (0.099) &   99.0 (0.128) &  98.2 (0.130)\\
&  & & $\sigma$ &-0.004 &  0.002 &   88.2 (0.137) &  87.6 (0.137) &   93.0 (0.163) &  93.2 (0.163) &   96.2 (0.215) &  97.0 (0.218)\\
  \cline{3-12}
    & & $(59^\star0,15)$ & $\gamma_0$ &-2.195 & 5.860 &   37.5 (3.637) &  9.7 (4.82) &   50.8 (4.333) &  14.9 (5.763) &   79.3 (5.695) &  33.6 (7.659)\\
  &  & & $\gamma_1$ &0.064 & 0.005 &   38.4 (0.108) &  10.2 (0.147) &   52.2 (0.128) &  15.9 (0.176) &   80.6 (0.168) &  35.5 (0.234)\\
&  & & $\sigma$ &0.054 & 0.005 &   88.7 (0.186) &  46.9 (0.271) &   96.6 (0.222) &  62.3 (0.325) &   99.8 (0.291) &  83.9 (0.437)\\
  \cline{3-12}
    & & $(15^\star(0,0,1,0))$ & $\gamma_0$ &0.019 & 0.777 &   87.3 (2.811) &  89.9 (3.000) &   93.8 (3.349) &  95.1 (3.562) &   97.8 (4.402) &  98.9 (4.649)\\
  &  & & $\gamma_1$ &-0.001 & 0.001 &   87.1 (0.082) &  89.8 (0.087) &   93.9 (0.098) &  94.9 (0.104) &   97.8 (0.129) &  99.0 (0.136)\\
&  & & $\sigma$ &-0.002 & 0.002 &   86.7 (0.131) &  88.4 (0.134) &   92.2 (0.157) &  93.0 (0.159) &   95.5 (0.206) &  96.6 (0.210)\\
 \hline
 \end{tabular}}
 \end{center}
 \label{t4}
 \end{table*}

\section{Conclusions and future research}\label{s8}

In this paper we have presented the Multiple Step-Stress model with Type-II and progressive Type-II censoring, useful to evaluate the reliability of highly reliable products in a fast and economical manner. Assuming lognormally distributed lifetimes, we have derived the maximum likelihood estimates of the parameters, the Fisher Information matrix for both types of censoring schemes. and the approximate and bootstrap confidence intervals. A Monte Carlo simulation study has been conducted in order to evaluate the performance of the proposed method. The results are shown by means of figures and tables to help analyze them.

It is shown that the method provides accurate and almost unbiased estimates, as well as well-performing confidence intervals. Better results are obtained when the number of steps in the model is large (as more precise results are obtained) and complete right censoring is avoided. As usual, the larger number of non-censored units, the more precise are the results.

As part of future work, we plan to apply the Multiple Step-Stress model with progressive Type-II censoring (we remind the reader that this type of censoring scheme is not common in the literature but it provides better results than the classical Type-II censoring) to the Generalized Gamma distribution, since many distributions commonly used in this area (such as the Lognormal distribution, the Exponential distribution, the Weibull distribution and the Gamma distribution) are special cases of the generalized gamma. This will permit having a general framework for several different distributions and determining which parametric model is appropriate for a given set of data, when it is not known in advance.

\appendices
\section{Notation}

\ap{$\gamma_0,\gamma_1$}{parameters of the link function}
\ap{$\Phi(\cdot)$}{standard normal CDF}
\ap{$\phi(\cdot)$}{standard normal PDF}
\ap{$\mu_i$}{location parameter of the lognormal distribution}
\ap{$\sigma$}{common scale parameter of the lognormal}
\ap{}{distribution}
\ap{$\bm\theta$}{vector of parameters}
\ap{$\tau_i$}{time at which stress level is changed from}
\ap{}{$x_i$ to $x_{i+1}$}
\ap{$A, E$}{parameters of the Arrhenius law}
\ap{$f_i$}{PDF for the lognormal distribution at stress level $i$}
\ap{$F_i$}{CDF for the lognormal distribution at stresslevel $i$}
\ap{${\cal F}$}{Fisher information matrix}
\ap{$g$}{PDF for the lognormal $k$-step stress model}
\ap{$G$}{CDF for the lognormal $k$-step stress model}
\ap{${\cal I}_{obs}$}{observed Fisher information matrix}
\ap{$k$}{Boltzmann's constant}
\ap{$L(\cdot)$}{likelihood function}
\ap{$l(\cdot)$}{log-likelihood function}
\ap{$m$}{number of stress levels}
\ap{$n$}{sample size}
\ap{$n_i$}{number of failures that occur at stress level $x_i$}
\ap{$r$}{number of non-censored units}
\ap{$R_k$}{number of censored units when the $k$-th  failure}
\ap{}{takes place}
\ap{$t_{i,j}$}{$j$-th ordered failure time of $n_i$ units at stress }
\ap{}{level $x_i$}
\ap{$T_i$}{lifetime of a test unit at stress level $x_i$}
\ap{$V$}{temperature in Kelvin degrees}
\ap{$x_i$}{stress level at step $i$}

\section*{Acknowledgment}
The author is indebted to Prof. N. Balakrishnan for his insightful comments. This work
was partly supported by the Spanish Government through project Project TRA2010-
17818 (research short stays program).

\bibliographystyle{IEEEtran}


\begin{thebibliography}{10}
\providecommand{\url}[1]{#1}
\csname url@samestyle\endcsname
\providecommand{\newblock}{\relax}
\providecommand{\bibinfo}[2]{#2}
\providecommand{\BIBentrySTDinterwordspacing}{\spaceskip=0pt\relax}
\providecommand{\BIBentryALTinterwordstretchfactor}{4}
\providecommand{\BIBentryALTinterwordspacing}{\spaceskip=\fontdimen2\font plus
\BIBentryALTinterwordstretchfactor\fontdimen3\font minus
  \fontdimen4\font\relax}
\providecommand{\BIBforeignlanguage}[2]{{%
\expandafter\ifx\csname l@#1\endcsname\relax
\typeout{** WARNING: IEEEtran.bst: No hyphenation pattern has been}%
\typeout{** loaded for the language `#1'. Using the pattern for}%
\typeout{** the default language instead.}%
\else
\language=\csname l@#1\endcsname
\fi
#2}}
\providecommand{\BIBdecl}{\relax}
\BIBdecl

\bibitem{LiuQ:2011}
X.~Liu and Q.~W. S., ``Modeling and planning of step-stress accelerated life
  tests with independent competing risks,'' vol.~60, no.~4, pp. 712--720, 2011.

\bibitem{XuBT:14}
A.~Xu, S.~Basu, and Y.~Tang, ``A full bayesian approach for masked data in
  step-stress accelerated life testing,'' vol.~63, no.~3, pp. 798--806, 2014.

\bibitem{Sedyakin:66}
N.~M. Sedyakin, ``On one physical principle in reliability theory (in
  {R}ussian),'' vol.~3, pp. 80--87, 1966.

\bibitem{Nelson:80}
W.~Nelson, ``Accelerated life testing: {S}tep-stress models and data
  analysis,'' vol.~29, pp. 103--108, 1980.

\bibitem{BagdonaviciusN:02}
V.~Bagdonavicius and M.~Nikulin, \emph{Accelerated Life Models: Modelling and
  Statistical Analysis}.\hskip 1em plus 0.5em minus 0.4em\relax Boca Raton,
  Florida: Chapman and Hall/CRC Press, 2002.

\bibitem{GangulyKM:15}
A.~Ganguly, D.~Kundu, and S.~Mitra, ``Bayesian analysis of a simple step-stress
  model under weibull lifetimes,'' vol.~64, no.~1, pp. 473--485, 2015.

\bibitem{Xiong:98}
C.~Xiong, ``Inference on a simple step-stress model with {T}ype-{II} censored
  exponential data,'' vol.~47, pp. 142--146, 1998.

\bibitem{BalakrishnanKNK:07}
N.~Balakrishnan, D.~Kundu, H.~NG, and N.~Kannan, ``Point andvinterval
  estimation for a simple step-stress model with {T}ype-{II} censoring,''
  vol.~39, pp. 35--47, 2007.

\bibitem{BalakrishnanX:07b}
N.~Balakrishnan and Q.~Xie, ``Exact inference for a simple step-stress model
  with {T}ype-{I} hybrid censored data from the exponential distribution,''
  vol. 137, pp. 3268--3290, 2007.

\bibitem{BalakrishnanX:07a}
------, ``Exact inference for a simple step-stress model with {T}ype-{II}
  hybrid censored data from the exponential distribution,'' vol. 137, pp.
  2543--2563, 2007.

\bibitem{BalakrishnanZX:09}
N.~Balakrishnan, L.~Zhang, and Q.~Xie, ``Inference for a simple step-stress
  model with {Type-I} censoring and lognormally distributed lifetimes,''
  vol.~38, pp. 1690--1709, 2009.

\bibitem{MillerN:83}
R.~Miller and W.~B. Nelson, ``Optimum simple step-stress plans for accelerated
  life testing,'' vol.~32, pp. 59--65, 1983.

\bibitem{Bai:89}
D.~Bai, M.~Kim, and S.~Lee, ``Optimum simple step-stress accelerated life test
  with censoring,'' vol.~38, pp. 528--532, 1989.

\bibitem{XieBD:07}
Q.~Xie, N.~Balakrishnan, and D.~Han, ``Advances in {M}athematical and
  {S}tatistical modeling,'' in \emph{Exact Inference and Optimal Censoring
  Scheme for a Simple Step-Stress Model Under Progressive {T}ype-{II}
  Censoring}, B.~C. Arnold, N.~Balakrishnan, J.~M. Sarabia, and R.~M\'inguez,
  Eds.\hskip 1em plus 0.5em minus 0.4em\relax Berlin: Birkhauser, 2009, pp.
  107--137.

\bibitem{KhamisH:98}
I.~H. Khamis and J.~J. Higgins, ``A new model for step-stress testing,''
  vol.~47, pp. 131--134, 1998.

\bibitem{KateriB:08}
M.~Kateri and N.~Balakrishnan, ``Inference for a simple step-stress model with
  {Type-II} censoring, and weibull distributed lifetimes,'' vol.~57, pp.
  616--626, 2008.

\bibitem{NgCB:04}
H.~Ng, P.~Chan, and B.~Balakrishnan, ``Optimal progressive censoring plan for
  the {W}eibull distribution,'' vol.~46, pp. 470--481, 2004.

\bibitem{LinC:12}
C.~Lin and C.~Chou, ``Statistical inference for a lognormal step-stress model
  with {T}ype-{I} censoring,'' vol.~61, pp. 361--377, 2012.

\bibitem{GuonoB:01}
E.~Gouno and N.~Balakrishnan, ``Step-stress accelerated life test,'' in
  \emph{Advances in Reliability}, ser. Handbook of Statistics, N.~Balakrishnan
  and C.~Rao, Eds.\hskip 1em plus 0.5em minus 0.4em\relax Amsterdam:
  North-Holland, 2001, vol.~20, pp. 623--639.

\bibitem{Nelson:05a}
W.~Nelson, ``A bibliography of accelerated test plans, {P}art {I}-{O}verview,''
  vol.~54, pp. 194--197, 2005.

\bibitem{Nelson:05b}
------, ``A bibliography of accelerated test plans, {P}art {II}-{R}eference,''
  vol.~54, pp. 370--373, 2005.

\bibitem{Balakrishnan:09}
N.~Balakrishnan, ``A synthesis of exact inferential results for exponential
  step-stress models and associated optimal accelerated life-tests,'' vol.~69,
  pp. 351--396, 2009.

\bibitem{Balakrishnan:07}
------, ``Progressive censoring methodology: an appraisal,'' vol.~16, no.~2,
  pp. 211--259, 2007.

\bibitem{IUPAC}
A.~D. McNaught and A.~W. (IUPAC), \emph{Compendium of Chemical Terminology, 2nd
  ed. (the "Gold Book")}.\hskip 1em plus 0.5em minus 0.4em\relax Oxford:
  Blackwell Scientific Publications, 1997.

\bibitem{EfronT:93}
B.~Efron and R.~Tibshirani, \emph{An introduction to the Bootstrap}.\hskip 1em
  plus 0.5em minus 0.4em\relax New York: Chapman and Hall, 1993.

\bibitem{BalakrishnanS:95}
N.~Balakrishnan and R.~A. Sandhu, ``A simple simulational algorithm for
  generating progressive type-{II} censored samples,'' vol.~49, no.~2, pp.
  229--230, 1995.

\end{thebibliography}

\end{document}